\documentclass[10pt]{article}
\usepackage{amsfonts}
\usepackage[leqno]{amsmath}
\usepackage{graphicx}
\usepackage{latexsym}
\usepackage{amsmath,amsfonts,amssymb,amsthm,mathrsfs,euscript,makeidx,color}
\usepackage{enumerate}
\allowdisplaybreaks[1]
\usepackage{multirow}

\def\sqr#1#2{{\vcenter{\vbox{\hrule height.#2pt
				\hbox{\vrule width.#2pt height#1pt \kern#1pt \vrule width.#2pt}
				\hrule height.#2pt}}}}
\def\signed #1{{\unskip\nobreak\hfil\penalty50
		\hskip2em\hbox{}\nobreak\hfil#1
		\parfillskip=0pt \finalhyphendemerits=0 \par}}
\def\endpf{\signed {$\sqr69$}}

\def\3n{\negthinspace \negthinspace \negthinspace }
\def\2n{\negthinspace \negthinspace }
\def\1n{\negthinspace }
\def\bel{\begin{equation}\label}

\def\dbE{\mathbb{E}}
\def\dbF{\mathbb{F}}

\def\dbH{\mathbb{H}}

\def\dbP{\mathbb{P}}

\def\dbR{\mathbb{R}}
\def\dbS{\mathbb{S}}

\def\dbX{\mathbb{X}}
\def\dbY{\mathbb{Y}}
\def\dbZ{\mathbb{Z}}

\def\sD{\mathscr{D}}

\def\sR{\mathscr{R}}

\def\sU{\mathscr{U}}

\def\ds{\displaystyle}

\def\ns{\noalign{\ss}}
\def\a{\alpha}
\def\b{\beta }
\def\g{\gamma}
\def\d{\delta}
\def\e{\varepsilon}
\def\z{\zeta}

\def\l{\lambda}

\def\si{\sigma}
\def\t{\tau}
\def\f{\varphi}
\def\th{\theta}

\def\i{\infty}
\def\G{\Gamma}

\def\L{\Lambda}

\def\F{\Phi}
\def\O{\Omega}
\def\cD{{\cal D}}

\def\cF{{\cal F}}

\def\cN{{\cal N}}

\def\cS{{\cal S}}

\def\cl{{\cal l}}
\def\BF{{\bf F}}

\def\BH{{\bf H}}

\def\BQ{{\bf Q}}
\def\BR{{\bf R}}
\def\BS{{\bf S}}

\def\G{\Gamma}

\def\L{\Lambda}

\def\F{\varPhi}
\def\O{\Omega}

\def\BBPsi{\boldsymbol\Psi}

\def\no{\noindent}

\def\ss{\smallskip}
\def\ms{\medskip}

\def\q{\quad}
\def\qq{\qquad}
\def\hb{\hbox}

\def\h1{\outline{$1$}}

\def\hh1{\outline{$1$}}
\def\hh2{\outline{$2$}}
\def\hh3{\outline{$3$}}
\def\hh4{\outline{$4$}}
\def\hh5{\outline{$5$}}
\def\hh6{\outline{$6$}}
\def\hh7{\outline{$7$}}
\def\hh8{\outline{$8$}}
\def\hh9{\outline{$9$}}
\def\hh0{\outline{$0$}}

\def\liminf{\mathop{\underline{\rm lim}}}

\def\Ra{\mathop{\Rightarrow}}

\def\lan{{\langle}}
\def\ran{{\rangle}}

\def\h{\widehat}
\def\wt{\widetilde}

\def\cd{\cdot}

\def\cl{\overline}

\def\les{\leqslant}
\def\ges{\geqslant}

\def\({\Big (}
\def\){\Big )}
\def\[{\Big[}
\def\]{\Big]}
\def\lan{\langle}
\def\ran{\rangle}

\def\bde{\begin{definition}\label}
	\def\ede{\end{definition}}
\def\bel{\begin{equation}\label}
		\def\ee{\end{equation}}
	\def\bt{\begin{theorem}\label}
		\def\et{\end{theorem}}
	\def\bc{\begin{corollary}\label}
		\def\ec{\end{corollary}}
	\def\bl{\begin{lemma}\label}
		\def\el{\end{lemma}}
	\def\bp{\begin{proposition}\label}
		\def\ep{\end{proposition}}
	\def\bex{\begin{example}\label}
		\def\ex{\end{example}}
	\def\bas{\begin{assumption}}
		\def\eas{\end{assumption}}
	\def\br{\begin{remark}\label}
		\def\er{\end{remark}}
	\def\ba{\begin{array}}
		\def\ea{\end{array}}
	\def\ed{\end{document}}

\def\rf{\eqref}

\def\square#1{\vbox{\hrule\hbox{\vrule height#1%
			\kern#1\vrule}\hrule}}
\def\rectangle#1#2{\vbox{\hrule\hbox{\vrule height#1%
			\kern#2\vrule}\hrule}}

\font\tenbb=msbm10 \font\sevenbb=msbm7 \font\fivebb=msbm5

\newfam\bbfam
\scriptscriptfont\bbfam=\fivebb \textfont\bbfam=\tenbb
\scriptfont\bbfam=\sevenbb

\newtheorem{theorem}{Theorem}[section]
\newtheorem{corollary}[theorem]{Corollary}

\newtheorem{lemma}[theorem]{Lemma}
\newtheorem{proposition}[theorem]{Proposition}

\theoremstyle{definition}
\newtheorem{definition}[theorem]{Definition}

\newtheorem{remark}[theorem]{Remark}

\newtheorem{example}[theorem]{Example}

\makeatletter

\@addtoreset{equation}{section}
\makeatother

\usepackage{xcolor}
\makeatletter

\input pdf-trans
\newbox\qbox
\def\usecolor#1{\csname\string\color@#1\endcsname\space}
\newcommand\bordercolor[1]{\colsplit{1}{#1}}
\newcommand\fillcolor[1]{\colsplit{0}{#1}}
\newcommand\outline[1]{\leavevmode%
	\def\maltext{#1}%
	\setbox\qbox=\hbox{\maltext}%
	\boxgs{Q q 2 Tr \thickness\space w \fillcol\space \bordercol\space}{}%
	\copy\qbox%
}
\makeatother
\newcommand\colsplit[2]{\colorlet{tmpcolor}{#2}\edef\tmp{\usecolor{tmpcolor}}%
	\def\tmpB{}\expandafter\colsplithelp\tmp\relax%
	\ifnum0=#1\relax\edef\fillcol{\tmpB}\else\edef\bordercol{\tmpC}\fi}
\def\colsplithelp#1#2 #3\relax{%
	\edef\tmpB{\tmpB#1#2 }%
	\ifnum `#1>`9\relax\def\tmpC{#3}\else\colsplithelp#3\relax\fi
}
\bordercolor{black}
\fillcolor{white}
\def\thickness{.3}

\begin{document}

\title{\bf Stochastic Optimal Linear Quadratic Controls with\\ A Recursive Cost Functional}

\author{Lin Li\thanks{Department of Mathematics, University of Central Florida, Orlando, FL 32816, USA (Email: {\tt lin.li@ucf.edu}).
                           }
~~~and~~~
Jiongmin Yong\thanks{Department of Mathematics, University of Central Florida, Orlando, FL 32816, USA
                    (Email: {\tt jiongmin.yong@ucf.edu}).
                    This author is supported by NSF grant DMS-2305475.}
}

\date{}
\maketitle
\centerline{({\it In the Memory of Professor Thomas I. Seidman})}
\ms

\no{\bf Abstract.} This paper is concerned with a stochastic linear quadratic (LQ, for short) control problem with a recursive cost functional. It involves BSDEs in $L^1$ whose well-posedness is a subtle issue. A suitable framework has been adopted so that the corresponding LQ problem is correctly formulated. Open-loop and closed-loop solvability of such an LQ problem have been investigated and characterized by the solvability of an FBSDE and that of Riccati differential equation.

\ms

\no{\bf Keywords.} Optimal linear quadratic control, recursive cost functional, forward-backward stochastic differential equation, Riccati equation, backward differential equation in $L^1$ space.

\ms

\no{\bf AMS 2020 Mathematics Subject Classification.} 93E20, 49N10, 60H10.

\section{Introduction.}\label{Sec:Intro}

Let $(\O,\cF,\dbP)$ be a complete probability space, with $\cN$ being the set of all $\dbP$-null sets, on which a standard Brownian motion is defined. Let $T>0$ be the given time horizon and let $t\in[0,T]$ be given. Define
$$\cF^t_s=\si\(\{W(\t')-W(\t)\bigm|t\les\t<\t'\les s\}\vee\cN\).$$
Denote $\dbF^t=\{\cF^t_s\}_{s\ges t}$. In what follows, we will denote $\dbE_s^t[\,\cd\,]=\dbE^t[\,\cd\,|\cF^t_s]$ (and thus $\dbE^t_t[\,\cd\,]=\dbE^t[\,\cd\,|\cF^t_t]=\dbE[\,\cd\,]$) Now, in $(\O,\cF,\dbF^t,\dbP)$, we consider the following controlled linear stochastic equation:
\bel{state}\left\{\2n\ba{ll}
\ds dX(s)=[A(s)X(s)+B(s)u(s)+b(s)]ds+[C(s)X(s)+D(s)u(s)+\si(s)]dW(s),\q s\in[t,T],\\
\ns\ds X(t)=x.\ea\right.\ee
For the above state equation, we introduce the following hypothesis:

\ms

{\bf(H1)} The maps $A,C:[0,T]\to\dbR^{n\times n}$, $B,D:[0,T]\to\dbR^{n\times m}$, and $b,\si:[0,T]\times\O\to\dbR^n$ satisfy the following:
$$\ba{ll}
\ns\ds A(\cd),C(\cd)\in L^\i(0,T;\dbR^{n\times n})\equiv\Big\{A:[0,T]\to\dbR^{n\times n}\bigm|A(\cd)\hb{ is measurable and bounded}\Big\},\\
\ns\ds B(\cd),D(\cd)\in L^\i(0,T;\dbR^{n\times m})\equiv\Big\{B:[0,T]\to\dbR^{n\times m}\bigm|B(\cd)\hb{ is measurable and bounded}\Big\},\\
\ns\ds b(\cd)\in L^1(0,T;\dbR^n)\equiv\Big\{b:[0,T]\to\dbR^n\bigm|b(\cd)\hb{ is measurable and }\int_t^T|b(s)|ds<\i\Big\},\\
\ns\ds\si(\cd)\in L^2(0,T;\dbR^n)\equiv\Big\{\si:[0,T]\to\dbR^n\bigm|\si(\cd)\hb{ is measurable and }\int_t^T|\si(s)|^2ds<\i\Big\},\ea$$
Note that in the above, we have assumed $b(\cd)$ and $\si(\cd)$ to be deterministic because they are required to be $\dbF^t$-progressively measurable for all $t\in[0,T]$ which leads to $b(\cd)$ and $\si(\cd)$ to be $\dbF^t$-progressively measurable for each $t\in[0,T]$, it implies that $b(t)$ and $\si(t)$ are $\si(\cN)$-measurable, for each $t\in[0,T]$. Thus, it is almost surely that $b(t)$ and $\si(t)$ are deterministic vectors for each $t\in[0,T]$. Hence, we assume that $b(\cd)$ and $\si(\cd)$ are deterministic, for convenience. Next, we introduce the set of admissible controls:
$$\ba{ll}
\ns\ds\sU^p[t,T]=\Big\{u:[t,T]\times\O\to\dbR^m\bigm|u(\cd)\hb{ is $\dbF$-progressively meaurble, }\dbE\(\int_t^T|u(s)|^2ds\)^{p\over2}<\i\Big\}\\
\ns\ds\qq\qq\equiv L^p_\dbF(\O;L^2(t,T;\dbR^m)),\qq p\ges1.\ea$$
Then, under (H1), for any initial pair $(t,x)\in[0,T)\times\dbR^n$ and $u(\cd)\in\sU^p[t,T]$, there exists a unique solution
$$\ba{ll}
\ns\ds X(\cd)=X(\cd\,;t,x,u(\cd))\in L^p_\dbF(\O;C[t,T];\dbR^n))\\
\ns\ds\equiv\Big\{X(\cd)\in L^p_\dbF(\O;L^1(t,T;\dbR^n))\bigm|X(\cd)\hb{ has continuous paths, }\dbE\[\sup_{s\in[t,T]}|X(s)|^p\]<\i\Big\},\ea$$
to the state equation \rf{state}. Moreover, the following holds:
\bel{|X|}\dbE\(\sup_{s\in[t,T]}|X(s)|^p\)\les K\[|x|^p+\(\int_t^T|b(s)|ds\)^p+\(\int_t^T|\si(s)|^2ds\)^{p\over2}+
\dbE\(\int_t^T|u(s)|^2ds\)^{p\over2}\],\ee
which implies $X(T)\in L^p_{\cF_T}(\O;\dbR^n)$. In what follows, we let $K>0$ to be a generic constant, which could be different from line to line. Further, we have the following stability estimate: If $u'(\cd)\in\sU^p[t,T]$ is another control with $X'(\cd)$ being the corresponding state process. Then
\bel{X-X}\dbE\(\sup_{s\in[t,T]}|X(s)-X'(s)|^p\)\les K\dbE\(\int_t^T|u(s)-u'(s)|^2ds\)^{p\over2}.\ee
To measure the performance of the controls, we introduce the following cost functional of the recursive form
\bel{cost}J(t,x;u(\cd))=Y(t),\ee
where $(Y(\cd),Z(\cd))$ is the adapted solution of the following stochastic backward differential equation (BSDE, for short):
\bel{BSDE1}\left\{\2n\ba{ll}
\ns\ds dY(s)=-[E(s)Y(s)+F(s)Z(s)+f(s)\ran]ds+Z(s)dW(s),\qq s\in[t,T],\\
\ns\ds Y(T)=\xi,\ea\right.\ee
where
\bel{f**}\ba{ll}
\ns\ds f(s)=\lan Q(s)X(s),X(s)\ran+2\lan S(s)X(s),u(s)\ran+\lan R(s)u(s),u(s)\ran+2\lan q(s),X(s)\ran+2\lan r(s),u(s)\ran,\\
\ns\ds\qq\qq\qq\qq\qq\qq\qq\qq\qq\qq\qq\qq\qq\qq\qq\qq s\in[t,T],\\
\ns\ds\xi=\lan GX(T),X(T)\ran+2\lan g,X(T)\ran,\ea\ee
for some deterministic bounded functions $E(\cd)$ and $F(\cd)$. Note that by assuming $Q(s)$, $S(s)$ and $R(s)$ to be deterministic and bounded, $q(s)$ and $r(s)$ to be deterministic square integrable processes, and $g\in\dbR^n$       (The argument for $q(\cd),r(\cd)$ and $g$ to be deterministic is similar to that for $b(\cd)$ and $\si(\cd)$)
$$\ba{ll}
\ns\ds\dbE\int_t^T|f(s)|ds\les K\dbE\int_t^T\(|X(s)|^2+|X(s)||u(s)|+|u(s)|^2+|X(s)|
+|u(s)|\)ds\\
\ns\ds\les K\dbE\int_t^T\(|X(s)|^2+|u(s)|^2+|q(s)|^2+|r(s)|^2\)ds\\
\ns\ds\les K\[|x|^2+\(\int_t^T|b(s)|ds\)^2+\int_t^T\(|\si(s)|^2+|q(s)|^2+|r(s)|^2\)ds
+\dbE\int_t^T|u(s)|^2ds\]<\i,\ea$$
and further assuming $\begin{pmatrix}Q(s)&S(s)^\top\\ S(s)&R(s)\end{pmatrix}$ to be uniformly positive definite, which leads to, for some $\d>0$,
$$\ba{ll}
\ns\ds\dbE\int_t^T|f(s)|ds\ges\d\dbE\int_t^T\(|X(s)|^2+|u(s)|^2\)ds-K\dbE\int_t^T
\(|X(s)||u(s)|+|q(s)||X(s)|+|r(s)|u(s)|\)ds\\
\ns\ds\ges{\d\over2}\dbE\int_t^T\(|X(s)|^2+|u(s)|^2\)ds-K\int_t^T\(|q(s)|^2
+|r(s)|^2\)ds\\
\ns\ds\ges{\d\over2}\dbE\int_t^T|u(s)|^2ds-K\int_t^T\(|q(s)|^2
+|r(s)|^2\)ds.\ea$$
Combining the above, we have
\bel{|f|}\ba{ll}
\ns\ds{\d\over2}\dbE\int_t^T|u(s)|^2ds-K\int_t^T\(|q(s)|^2
+|r(s)|^2\)ds\les\dbE\int_t^T|f(s)|ds\\
\ns\ds\les K\[|x|^2+\(\int_t^T|b(s)|ds\)^2+\int_t^T\(|\si(s)|^2+|q(s)|^2+|r(s)|^2\)ds
+\dbE\int_t^T|u(s)|^2ds\].\ea\ee
Thus, essentially, $\|f(\cd)\|_{L_
\dbF^1(\O;L^1(t,T:\dbR^n))}$ and $\|u(\cd)\|_{\sU^2[t,T]}$ are of the same orders. Next, by assuming square integrability of $g$, we get
$$\ba{ll}
\ns\ds\dbE|\xi|\les K\dbE\(|X(T)|^2+|g|^2\)\\
\ns\ds\les K\dbE\[|x|^2+\(\int_t^T|b(s)|ds\)^2
+\int_t^T|\si(s)|^2ds+\int_t^T|u(s)|^2ds+|g|^2\]<\i.\ea$$
Further, let $G$ to be positive definite and
$$B(s)=I,\q D(s)=0,\q b(s)=\si(s)=0,\qq s\in[t,T],\qq x=0.$$
For $k\ges1$, we let
$$u_k(s)={1\over(T-t_k)^{1\over2}}{\bf1}_{[t_k,T]}(s)u_0,\qq s\in[t,T],$$
for some fixed $u_0\in\dbR^m\setminus\{0\}$. Then
$$X_k(T)=\int_{t_k}^T\Psi(\t,t_k)u_k(\t)d\t,\qq\dbE\int_t^T|u_k(s)|^2ds=|u_0|^2.$$
with $\Psi(\cd\,,\cd)$ solving the following
$$\left\{\2n\ba{ll}
\ns\ds d\Psi(t,t_k)=A(s)\Psi(s,t_k)ds+C(s)\Psi(s,t_k)dW(s),\\
\ns\ds\Psi(t_k,t_k)=I.\ea\right.$$
Hence,
$$\dbE|\xi|\ges\dbE\(\d|X_k(T)|^2-|g|^2\)\ges{\d\over T-t_k}\dbE\Big|\int_{t_k}^T\Psi(\t,t_k)u_0
d\t\Big|^2=\d\dbE\int_t^T|u_k(s)|^2ds+o(1).$$
Combing the above, in general, we should have
\bel{|xi|}\d\dbE\int_t^T\3n|u_k(s)|^2ds+o(1)\les\dbE|\xi|\les K\dbE\[|x|^2+\(\int_t^T\3n|b(s)|ds\)^2
+\int_t^T\3n|\si(s)|^2ds+\int_t^T\3n|u_k(s)|^2ds+|g|^2\].\ee
Thus, in general, $\|\xi\|_{L^1_{\cF_T}(\O:\dbR^n)}$ and $\|u(\cd)\|_{\sU^2[t,T]}$ are of the same orders. The above shows that, in general, \rf{BSDE1} is a linear nonhomogeneous BSDE with the nonhomogeneous term being an $L^1$ process and the terminal term being an $L^1$ random variable. Thus, the well-posedness of such a BSDE is very subtle (see \cite{Briand-Hu 2005, Fan-Liu 2010, Fan 2016, Klimsaik-Rzymowski 2024} among many papers). We will see that for $u(\cd)\in\sD(J(t,x;\cd))$, a dense subset in $\sU^2[t,T]$, the BDSE \rf{BSDE1} is well-psosed, and the cost functional \rf{cost} is well-defined which is called a {\it recursive functional}. The recursive functional was firstly introduced by Duffie--Epstein in 1992 (\cite{Duffie-Epstein 1992a,Duffie-Epstein 1992b}) (which was called differential utility). People found that this was nothing but some special cases of BSDEs. Peng then firstly studied some general nonlinear optimal control problems with recursive cost functional in 1992 (\cite{Peng 1992}). A little later, Peng introduced the so-called $g$-expectation notion (a generalized expectation, which could nonlinear, but is time-consistent) in 1997 from general theory of BSDEs (\cite{Peng 1997}). The notion of differential utility was extended by Quenez--Lazrak \cite{Lazrak-Quenez 2003}, and Lazrak \cite{Lazrak 2004} in the beginning of 2000s to a general BSDE framework. In 2008, Yong studied LQ problems under (nonlinear) $g$-expectation introduced in \cite{Peng 1997} with the framework of $L^{p\over2}$ space ($p>2$) \cite{Yong 2008}). Thus, the linear structure was partially ruined (since the expectation was nonlinear) and moreover, it could not cover the case of $p=2$. Because of that, the Riccati type differential equation was derived only for a very special case in \cite{Yong 2008}. In the current paper, with a recursive (in $L^1$, and linear) cost functional, defined through the idea of \cite{Peng 1997}, we are able to pose the following optimal control problem for which all the linear-quadratical structure is kept.

\ms

\bf Problem (LQR). \rm For given initial pair $(t,x)\in[0,T]\times\dbR^n$, find a control $\bar u(\cd)\in\sD(J(t,x;\cd))\subseteq\sU^2[t,T]$ such that
\bel{bar u}J(t,x;\bar u(\cd))=\inf_{u(\cd)\in\sD(J(t,x;\cd))}J(t.x;u(\cd))=V(t,x).\ee

\ms

Note that when $E(s),F(s)\equiv0$, the cost reads
\bel{}\ba{ll}
\ns\ds J(t,x;u(\cd))=\dbE\[\lan GX(T),X(T)\ran+2\lan g,X(T)\ran+\int_t^T\(\lan Q(s)X(s),X(s)\ran+2\lan S(s)X(s),u(s)\ran\\
\ns\ds\qq\qq\qq\qq\qq\qq\qq\qq+\lan R(s)u(s),u(s)\ran+2\lan q(s),X(s)\ran+2\lan r(s),u(s)\ran\)ds\].\ea\ee
In this case, we actually have $\sD(J(t,x;\cd))=\sU[t,T]^2$. Then Problem (LQR) is well-posed for $p=2$ and the above gives a classical LQ cost functional. Therefore, our results will cover the classical LQ theory presented in, say, \cite{Sun-Yong 2020}. Now, we highlight the main features of the current paper.

\ms

$\bullet$ Developping the BSDE in $L^1$ space (which is an extension of \cite{Peng 1997}), we correctly formulate the stochastic Problem (LQR) with $L^1$ recursive cost functional.

\ms

$\bullet$ Revealing a relation of Problem (LQR) and the classical LQ problem. It turns out that the former under a transformation leads to a classical LQ problem. But the latter is not a simple extension of the former, as far as the cost functionals (together with their domains) are concerned.

\ms

$\bullet$ Introduce the open-loop solvability of Problem (LQR) and its equivalent conditions, which is the solvability of an FBSDE, satisfying the maximum condition. This will be achieved by the idea of \cite{Mou-Yong 2006}. We will find the Fr\'echet derivative and the Hessian of the recursive cost functional.

\ms

$\bullet$ Introduce the closed-loop solvability of Problem (LQR) and its equivalent conditions, which is the solvability of a Riccati differential equation. The basic approach is the completing the squares.

\ms

The rest of this paper is organized as follows, Section 2 is to present some relevant results of BSDEs in $L^1$, together with the formulation of our LQ problem, and open-loop, and closed-loop solvability of the problem. By making some transformation, it seems that a classical LQ is derived. But, this LQ problem is not simple extension of our Problem (LQR). A concrete relation will be revealed in Section 3. In Section 4, we present the equivalence between open-loop solvability and that of an FBSDE. Closed-loop solvability of the Problem (LQ) turns out to be equivalent to that of Riccati equation. This will be presented in Section 5. Finally, some conclusion remarks are collected in Section 6.

\section{BSDE in $L^1$ and the Formulation of Problem (LQR).}

We first introduce the following definitions.

\bde{} \rm Let $\b\in(0,1)$. Define:
$$\ba{ll}
\ns\ds\cS^\b[t,T]=\Big\{Y:[t,T]\times\O\to\dbR\bigm|Y(\cd)\hb{ is $\dbF$-progressively measurable,}\\
\ns\ds\qq\qq\qq\qq\qq\qq\qq\|Y(\cd)\|_{\cS^\b[t,T]}\equiv\[\dbE\(\sup_{s\in[t,T]}
|Y(s)|^\b\)\]^{1\over\b}<\i\Big\},\\
\ns\ds M^\b[t,T]=\Big\{Z:[t,T]\times\O\to\dbR\bigm|Z(\cd)\hb{ is $\dbF$-progressively measurable,}\\
\ns\ds\qq\qq\qq\qq\qq\qq\qq\|Z\|_{M^\b}=\[\dbE\(\int_t^T|Z(s)|^2ds\)^{\b\over2}
\]^{1\over\b}<\i\Big\}.\ea$$

\ede

Clearly, $\cS^\b[t,T]$ and $M^\b[t,T]$ are complete linear metric spaces.

\bde{} \rm Process $Y:[t,T]\times\O\to\dbR$ is of class (D) if it is $\dbF$-adapted, and
$$\{Y(\t)\bigm|\t\hb{ is a stopping time valued in $[t,T]$ }\}$$
is uniformly integrable.

\ede

Consider the following general BSDE:
\bel{BSDE-g}\left\{\2n\ba{ll}
\ds dY(s)=-g(s,Y(s),Z(s))ds+Z(s)dW(s),\qq s\in[t,T],\\
\ns\ds Y(T)=\xi.\ea\right.\ee
We assume the following.

\ms

{\bf(H2)} Let
\bel{E[|xi|]<i}\dbE\[|\xi|+\int_0^T|g(s,0,0)|ds\]<\i,\ee
and for some constant $K>0$ and $\a\in(0,1)$,
\bel{g-g}|g(s,y_1,z_1)-g(s,y_2,z_2)|\les K\(|y_1-y_2|+|z_1-z_2|^\a\),\qq y_1,y_2,z_1,z_2\in\dbR,\ee
for some $\a\in(0,1)$.

\ms

According to \cite{Fan-Liu 2010} (and \cite{Briand-Hu 2005} p.16, see also the introduction of \cite{Fan 2016}, and \cite{Klimsaik-Rzymowski 2024}), we have the following result.\footnote{The authors would like to thank Professor Shengjun Fan (China University of Mining and Technology), and Professor Jianfeng Zhang (University of Southern California) for clarification of the results for general nonlinear BSDEs in $L^1$.}

\bt{Th2.3} \sl Let {\rm(H2)} hold. Then BSDE admits a unique adapted solution $(Y(\cd),Z(\cd))$, with $Y(\cd)$ in the class (D) and such that
\bel{2.4}\dbE\(\sup_{s\in[0,T]}|Y(s)|^\b\)^{1\over\b}+\dbE\(\int_t^T|Z(s)|^2
ds\)^{\b\over2}<\i,\qq\forall\b\in(0,1).\ee

\et

Note that the uniqueness follows the following comparison: If $\xi'\ges\xi$ are two proper random variables, and let $(Y'(\cd),Z'(\cd))$ and $(Y(\cd),Z(\cd))$ are corresponding adapted solutions, then
$$Y'(s)\ges Y(s),\qq s\in[t,T].$$
This comparison gives the uniqueness of the adapted solutions. The above results can be used to take care of the case where $F(s)\equiv0$ (see \rf{BSDE1}). However, this will not cover the case $F(s)\not\equiv0$. Hence, in what follows, we are going to adopting the notion of generalized solution from \cite{Peng 1997}, so that the case of $F(s)\not\equiv0$ can also be taken care.

\ms

In what follows, for any $p_1,p_2\ges1$ and Euclidean space $\dbH$, we denote
\bel{dbH}\ba{ll}
\ns\ds L^{p_1}_\dbF(\O;L^{p_2}(t,T;\dbH))=\Big\{\f:[t,T]\times\O\to\dbH\bigm|\f(\cd)\hb{ is $\dbF$-progressively measurable, }\\
\ns\ds\qq\qq\qq\qq\qq\qq\qq\qq\qq\qq\dbE\(\int_t^T|\f(s)|^{p_2}ds\)^{p_1\over p_2}<\i\Big\}.\ea\ee
\bel{dbH*}L^{p_1}_{\cF_T}(\O;\dbH)=\Big\{\xi:\O\to\dbH\bigm|\xi\hb{ is $\cF_T$-measurable, }\dbE|\xi|^{p_1}<\i\Big\}.\ee
We define
\bel{L^1+}L^{p_1+}_\dbF(\O;L^{p_2+}(t,T;\dbR))=\bigcup_{p_1'>p_1,p_2'>p_2}L^{p_1'}
(\O;L^{p_2'}(t,T;\dbR)),
\qq L^{p_1+}_{\cF_T}(\O;\dbR))=\bigcup_{p_1'>p_1}L^{p_1'}_{\cF_T}(\O;\dbR),\ee
and
\bel{sU^2+}\sU^{2+}[t,T]=\bigcup_{p>2}\sU^p[t,T].\ee
We consider BSDE \rf{BSDE1} with the general $f(\cd)$ and $\xi$ (they do not have to be forms \rf{f**}). Thus, \rf{BSDE1} is a linear BSDE with
\bel{f1}f(\cd)\in L^1_\dbF(\O;L^1(t,T;\dbR)),\qq \xi\in L^1_{\cF_T}(\O;\dbR).\ee
It is known that the well-posedness of BSDE in $L^1$ is very tricky. In the current paper, we are going to use the idea of \cite{Peng 1997} to approach the problem. The following lemma will play an essential role below.

\bl{} \sl Let $E(\cd),F(\cd)\in L^\i(t,T;\dbR)$. Let $\F(\cd\,,t)$ solve the following BSDE:
\bel{dF}\left\{\2n\ba{ll}
\ds d\F(s,t)=E(s)\F(s,t)ds+F(s)\F(s,t)dW(s),\qq0\les t\les s\les T,\\
\ns\ds\F(t,t)=1.\ea\right.\ee
Then
\bel{F}\F(s,t)=e^{\int_t^sE(\th)d\th}\[e^{\int_t^sF(\th)dW(\th)-{1\over2}\int_t^s
F(\th)^2d\th}\],\qq 0\les t,s\les T,\ee
and for any
\bel{f}f(\cd)\in L^{1+}_\dbF(\O;L^{1+}(t,T;\dbR)),\qq\xi\in L^{1+}_{\cF_T}(\O;\dbR),\ee
the unique adapted solution $(Y(\cd),Z(\cd))$ of BSDE \rf{BSDE1} admits the following representation:
\bel{Y}Y(s)=\dbE_s^t\(\F(T,s)\xi+\int_s^T
\F(s,\t)f(\t)d\t\),\qq s\in[t,T].\ee
Conversely, if \rf{Y} holds for some $\xi$ and $f(\cd)$ satisfying \rf{f}, then for some $Z(\cd)\in L^{1+}_\dbF(\O;L^{2+}(t,T;\dbR^n))$ such that $(Y(\cd),Z(\cd))$ is the adapted solution to BSDE \rf{BSDE1}.

\el

\it Proof. \rm  BDSDE \rf{dF} is a classical one, it is well-posed. Thus we may let $\F(\cd\,,t)$ be the solution of \rf{dF}. Set
$$\h\F(s,t)=e^{-\int_t^sE(\th)d\th}\F(s,t),\qq0\les t\les s\les T.$$
Then, \rf{dF} is equivalent to the following:
\bel{dhF}\left\{\2n\ba{ll}
\ds d\h\F(s,t)=F(s)\h\F(s,t)dW(s),\qq0\les t\les s\les T,\\
\ns\ds\h\F(t,t)=1.\ea\right.\ee
Thus,
$$d\(\ln\h\F(s,t)\)={1\over\h\F(s,t)}\[F(s)\h\F(s,t)dW(s)\]-{1\over2}
{1\over\h\F(s,t)^2}F(s)^2
\h\F(s,t)^2ds=F(s)dW(s)-{1\over2}F(s)^2ds.$$
Consequently, \rf{F} follows and the following group type property holds:
$$\h\F(s,s')\h\F(s',s'')=\h\F(s,s''),\qq0\les t\les s,s',s''\les T.$$
Now, for some $p>2$, we take
\bel{f-p}f(\cd)\in L^{p\over2}_\dbF(\O;L^{p\over2}(t,T;\dbR)),\qq\xi\in L^{p\over2}_{\cF_T}(\O;\dbR),\ee
and consider the corresponding BSDE \rf{BSDE1}.
Clearly, this BSDE is equivalent to the following:
\bel{BSDE-hY}\left\{\2n\ba{ll}
\ds d\h Y(s)=-[F(s)\h Z(s)+\h f(s)]ds+\h Z(s)dW(s),\qq s\in[t,T],\\
\ns\ds\h Y(T)=\h\xi,\ea\right.\ee
with
$$\h f(s)=e^{-\int_t^sE(\th)d\th}f(s),\q\h\xi=e^{-\int_t^TE(\th)d\th}\xi.$$
Since BSDE \rf{BSDE-hY} is well-posed, we may let $(\h Y(\cd),\h Z(\cd))$ be its adapted solution. Applying It\^o's formula to the map $s\mapsto\h\F(s,t)\h Y(s)$, we get
$$\ba{ll}
\ns\ds d\(\h\F(s,t)\h Y(s)\)=\(F(s)\h\F(s,t)dW(s)\)\h Y(s)+\h\F(s,t)
\[-\(F(s)\h Z(s,t)+\h f(s)\)ds+\h Z(s,t)dW(s)\]\\
\ns\ds\qq\qq\qq\qq\qq\qq\qq\qq\qq+F(s)\h\F(s,t)\h Z(s,t)ds\\
\ns\ds\qq\qq\qq=-\h\F(s,t)\h f(s)ds+\h\F(s,t)\(\h Y(s,t)F(s)+\h Z(s,t)\)dW(s),\qq s\in[t,T].\ea$$
Thus,
$$\ba{ll}
\ns\ds\dbE_s^t\[\h\F(T,t)\h\xi\]-\h\F(s,t)\h Y(s)=-\dbE_s^t\int_s^T
\h\F(\t,t)\h f(\t)d\t\ea$$
Consequently,
\bel{hY}\ba{ll}
\ns\ds\h Y(s)=\h\F(s,t)^{-1}\dbE_s^t\(\big[\h\F(T,t)\xi\big]+
\int_s^T\h\F(\t,t)\h f(\t)d\t\)\\
\ns\ds\qq=\dbE_s^t\(\h\F(T,s)\h\xi+
\int_s^T\h\F(\t,s)\h f(\t)d\t\).\ea\ee
Now, set
$$Y(s)=e^{\int_t^sE(\th)d\th}\h Y(s),\q Z(s)=e^{-\int_t^sE(\th)d\th}\h Z(s).$$ %
Multiplying \rf{hY} by $e^{\int_t^sE(\th)d\th}$, we have
$$\ba{ll}
\ns\ds Y(s)=e^{\int_t^sE(\th)d\th}\dbE^t_s\(\h\F(T,s)\h\xi+\int_s^T\h\F(\t,s)\h\F(
\h f(\t)d\t\)\\
\ns\ds\qq=\dbE^t_s\(e^{\int_t^sE(\th)d\th}\big[e^{\int_s^TE(\th)d\th}\F(T,s)\big]
e^{-\int_t^TE(\th)d\th}\xi+\int_s^Te^{\int_t^sE(\th)d\th}\big[e^{\int_s^\t E(\th)d\th}\F(\t,s)\big]
e^{-\int_t^\t E(\th)d\th}f(\t)d\t\)\\
\ns\ds\qq=\dbE_s^t\(\F(T,s)\xi+\int_s^T\F(\t,s)f(\t)d\t\).\ea$$
This proves our lemma. \endpf

\ms

When $E(s)=F(s)=0$ for all $s\in[t,T]$, the above representation \rf{Y} is called the $g$-expectation introduced in \cite{Peng 1997}; for the case that $E(s)=0$ for all $s\in[t,T]$, the representation reads
\bel{Y*}Y(s)=\dbE_s^t\[\F(T,s)\xi+\int_s^T\F(\t,s)f(\t)d\t\],\qq s\in[t,T].\ee
This is different from that in \cite{Peng 1997} (The function $\F(\t,s)$ was missing in the integrand on the right which could be a typo, we believe.).
Next result is basically quoted from \cite{Yong 2006}.

\bl{E[F(T)} \sl Let $F(\cd)\in L^\i(t,T;\dbR)$. Then for any $p\ges0$,
\bel{E[pF]}\dbE\Big|\dbE^t_s\[e^{\int_t^TF(\th)dW(\th)}\]\Big|^p\les e^{{p^2\over2}\int_s^TF(\th)^2d\th}.\ee

\el

\it Proof. \rm By \cite{Yong 2006} p.123, we have (with $p'\ges0$)
$$\ba{ll}
\ns\ds\dbE\Big|\dbE^t_s\[e^{\int_t^TF(\th)dW(\th)}\]\Big|^p=\sum_{k=0}^\i{p^k\over k!}\dbE\[\int_s^TF(\th)dW(\th)\]^k=\sum_{k=0}^\i
{p^{2k}\over(2k)!}{(2k)!\over2^kk!}\(\int_s^T|F(\th)|^2d\th\)^k\\
\ns\ds=\sum_{k=0}^\i{1\over k!}\({p^2\over2}\int_s^T|F(\th)|^2d\th\)^k=
e^{{p^2\over2}\int_s^T|F(\th)|^2d\th}.\ea$$
Then, our conclusion follows. \endpf

\ms

Consequently, for $p>2$, and for $\xi\in L^{p\over2}_{\cF_T}(\O;\dbR)$,
\bel{Fxi}\ba{ll}
\ns\ds\dbE\Big|\dbE_s^t\[e^{\int_t^TF(\th)dW(\th)}\xi\]\Big|\les
\(\dbE\[e^{\int_t^TF(\th)dW(\th)}\]^{p\over p-2}\)^{p-2\over p}\(\dbE|\xi|^{p\over2}\)^{2\over p}\\
\ns\ds\les\Big\{\dbE\[e^{{p^2\over2(p-2)^2}\int_s^TF(\th)^2d\th}\]
\Big\}^{p-2\over p}\(\dbE|\xi|^{p\over2}\)^{2\over p}\les K\(\dbE|\xi|^{p\over2}\)^{2\over p}.\ea\ee
Likewise, for $f(\cd)\in L^{p\over2}_\dbF(\O;L^{p\over2}(t,T;\dbR))$ (with $p>2$),
\bel{Ff}\ba{ll}
\ns\ds\dbE\Big|\dbE_s^t\int_s^Te^{\int_t^\t F(\th)dW(\th)}f(\t)d\t\Big|\les\(\dbE\int_s^T\[e^{\int_t^\t F(\th)dW(\th)}\]^{p\over p-2}d\t\)^{p-2\over p}\(\dbE\int_s^T|f(\t)|^{p\over2}d\t\)^{2\over p}\\
\ns\ds\les\(\int_s^Te^{{p^2\over2(p-2)^2}\int_s^\t F(\th)^2d\th}d\t\)^{p-2\over p}\(\dbE\int_s^T|f(\t)|^{p\over2}d\t\)^{2\over p}\les K\(\dbE\int_s^T|f(\t)|^{p\over2}d\t\)^{2\over p}.\ea\ee
Now, since $e^{\int_t^\t F(\th)dW(\th)}$ is unbounded, \rf{Y} might not be well-defined for all $f(\cd)$ and $\xi$ satisfying \rf{f1}. We thus introduce the following subsets.
$$\ba{ll}
\ns\ds\h L_\dbF^1(\O;L^1(t,T;\dbR))=\Big\{f(\cd)\in L^1_\dbF(\O;L^1(t,T;\dbR))\bigm|f^+(\cd),\hb{ or }f^-(\cd)\in L_\dbF^{1+}(\O;L^{1+}(t,T;\dbR))\Big\}\\
\ns\ds\qq\qq\qq\qq\qq\qq\subsetneqq L^1_\dbF(\O;L^1(t,T;\dbR)),\\
\ns\ds\h L_{\cF_T}^1(\O;\dbR)=\Big\{\xi\in L^1_{\cF_T}(\O;\dbR)\bigm|\xi^+,\hb{ or }\xi^-\in L^{1+}_{\cF_T}(\O;\dbR)\Big\}
\subsetneqq L^1_{\cF_T}(\O;\dbR),\ea$$
with $a^+=a\vee0$ and $a^-=(-a)\vee0$ (thus $a=a^+-a^-$). We point out that the above are not linear spaces. This can be seen from the following
$$\xi\in L^1_{\cF_T}(\O;\dbR)\setminus\h L^1_{\cF_T}(\O;\dbR)\q\Ra\q\xi^+,~\xi^-\in\h L^1_{\cF_T}(\O;\dbR).$$
But, one has
$$\xi=\xi^+-\xi^-\not\in\h L^1_{\cF_T}(\O;\dbR).$$
Note that in \cite{Peng 1997}, only non-negative process $f(\cd)$ and non-negative random variable $\xi$ are involved. Now, for any $f(\cd)\in\h L^1_\dbF(\O;L^1(t,T;\dbR))$ and $\xi\in\h L^1_{\cF_T}(\O;\dbR)$ with
\bel{f^-xi^-}f^-(\cd)\in L_\dbF^{1+}(\O;L^{1+}(t,T;\dbR)),\qq\xi^-\in
L^{1+}_{\cF_T}(\O;\dbR),\ee
we consider BSDE \rf{BSDE1}. Define
\bel{f_k}f^+_k(\cd)=f^+(\cd)\land k,\q\xi^+_k=\xi^+\land k,\qq k>0,\ee
and consider BSDEs:
\bel{BSDE4}\left\{\2n\ba{ll}
\ns\ds dY^+_k(s)=-[E(s)Y^+_k(s)+F(s)Z^+_k(s)+f_k^+(s)]ds+Z_k^+(s)dW(s),\qq s\in[t,T],\\
\ns\ds Y^+_k(T)=\xi^+_k,\ea\right.\ee
and by \rf{Fxi}--\rf{Ff}, the following BSDE is well-posed:
\bel{BSDE5}\left\{\2n\ba{ll}
\ns\ds dY^-(s)=-[E(s)Y^-(s)+F(s)Z^-(s)+f^-(s)]ds+Z^-(s)dW(s),\qq s\in[t,T],\\
\ns\ds Y^-(T)=\xi^-.\ea\right.\ee
We have the following representations:
$$Y^+_k(s)=\dbE_s^t\(\F(T,s)\xi^+_k+\int_s^T\F(\t,s)
f^+_k(\t)d\t\),\q s\in[t,T],$$
and
$$Y^-(s)=\dbE_s^t\(\F(T,s)\xi^-+\int_s^T
\F(\t,s)f^-(\t)d\t\),\q s\in[t,T].$$
The sequences $f^+_k(\cd)$ and $\xi^+_k$ are nondecreasing. Thus,
$$\lim_{k\to\i}f^+_k(\cd)=f^+(\cd),\qq\lim_{k\to\i}\xi^+_k=\xi^+.$$
By the comparison of BSDEs, we have
$$Y_k^+(s)\nearrow Y^+(s)=\dbE_s^t\(\F(T,s)\xi^++\int_s^T\F(\t,s)
f^+(\t)d\t\),\q s\in[t,T],$$
which might not be well-defined on $L^1_\dbF(\O;L^1(t,T;\dbR))\times L^1_{\cF_T}(\O;\dbR)$. On the other hand, by \rf{Fxi}--\rf{Ff},
$$Y^-(s)=\dbE_s\(\F(T,s)\xi^-+\int_s^T\F(\t,s)
f^-(\t)d\t\),\q s\in[t,T],$$
is well-defined. Hence, for any $f(\cd)\in\h L^1_\dbF(\O;L^1(t,T;\dbR))$ with $f^-(\cd)\in L_\dbF^{1+}(\O;L^{1+}(t,T;\dbR))$ and $\xi\in\h L^1_{\cF_T}(\O;\dbR)$ with $\xi^-\in L^{1+}_{\cF_T}(\O;\dbR)$, the {\it generalized solution} of the \rf{BSDE1} is defined as follows:
\bel{Y(s)}\ba{ll}
\ns\ds Y(s)=Y^+(s)-Y^-(s)\\
\ns\ds=\dbE_s\[\F(T,s)\(\xi^+-\xi^-\)
+\int_s^T\F(\t,s)
\(f^+(\t)-f^-(\t)\)d\t\],\q s\in[t,T],\ea\ee
Note that if \rf{f^-xi^-} holds, it is possible that
$$f^+(\cd)\notin L_\dbF^{1+}(\O;L^{1+}(t,T;\dbR)),\q\hb{or}\q\xi^+\notin
L^{1+}_{\cF_T}(\O;\dbR).$$
Thus, there are other three cases, for which the generalized solutions can be similarly defined:

\ms

$\bullet$ $f(\cd)\in\h L^1_\dbF(\O;L^1(t,T;\dbR))$ and $\xi\in\h L^1_{\cF_T}(\O;\dbR)$ with
\bel{f^-xi^+}f^-(\cd)\in L_\dbF^{1+}(\O;L^{1+}(t,T;\dbR)),
\qq\xi^+\in L^{1+}_{\cF_T}(\O;\dbR).\ee

\ms

$\bullet$ $f(\cd)\in\h L^1_\dbF(\O;L^1(t,T;\dbR))$ and $\xi\in\h L^1_{\cF_T}(\O;\dbR)$ with
\bel{f^+xi^-}f^+(\cd)\in L_\dbF^{1+}(\O;L^{1+}(t,T;\dbR)),\qq\xi^-\in L^{1+}_{\cF_T}(\O;\dbR).\ee

\ms

$\bullet$ $f(\cd)\in\h L^1_\dbF(\O;L^1(t,T;\dbR))$ and $\xi\in\h L^1_{\cF_T}(\O;\dbR)$ with
\bel{f^+xi^+}f^+(\cd)\in L_\dbF^{1+}(\O;L^{1+}(t,T;\dbR)),\qq\xi^+\in L^{1+}_{\cF_T}(\O;\dbR).\ee

\ms

Hence, the BSDE \rf{BSDE1} is solvable for $(f(\cd),\xi)\in\h L_\dbF^1(\O;L^1(t,T;\dbR))\times\h L^1_{\cF_T}(\O;\dbR)$. Consequently, we must have
$$\ba{ll}
\ns\ds L_\dbF^{1+}(\O;L^{1+}(t,T;\dbR))\times L^{1+}_{\cF_T}(\O;\dbR)\subseteq \h L_\dbF^1(\O;L^1(t,T;\dbR))\times\h L^1_{\cF_T}(\O;\dbR)\\
\ns\ds\qq\qq\qq\qq\qq\subseteq L_\dbF^1(\O;L^1(t,T;\dbR))\times L^1_{\cF_T}(\O;\dbR).\ea$$

\ms

In what follows, we let $\dbS^n$, $\dbS^n_+$ and $\dbS^n_{++}$ be the sets of all $(n\times n)$ symmetric, positive semi-definite, positive definite matrices, respectively. Now, we introduce the following assumption, relevant to Problem (LQR).

\ms

{\bf(H3)} Let the following hold:
$$\ba{ll}
\ns\ds G\in\dbS^n_+,\q g\in\dbR^n,\\
\ns\ds E(\cd),F(\cd)\in L^\i(0,T;\dbR),\q Q(\cd)\in L^\i(0,T;\dbS^n_+),\q S(\cd)\in L^\i(0,T;\dbR^{n\times m}),\q R(\cd)\in L^\i(0,T;\dbS^m_{++}),\\
\ns\ds q(\cd)\in L^2(0,T;\dbR^n),\q r(\cd)\in L^2(0,T;\dbR^m).\ea$$
Under (H1) and (H3), the state equation \rf{state} is well-posed, and \rf{f**} is defined. For $J(t,x;\cd)$ to be defined, we should take $u(\cd)$ from the set
\bel{sD}\sD(J(t,x;\cd))=\Big\{u(\cd)\in\sU^2[t,T]\bigm|J(t,x;u(\cd))\hb{ is well-defined and finite }\Big\},\ee
which is a linear space. Further, from Lemma \ref{E[F(T)} and \rf{Fxi}--\rf{Ff}, we see that (see \rf{sU^2+})
$$\sU^{2+}[t,T]\subseteq\Big\{u(\cd)\in\sU^2[t,T]\bigm|f(\cd)\in\h L_\dbF^1(\O;L^1(t,T;\dbR)),~\xi\in\h L^1_{\cF_T}(\O;\dbR)\hb{ with \rf{f**}}\Big\}\subseteq\sU^2[t,T].$$
In particular, for $f(\cd)$ and $\xi$ of form \rf{f**},
$$J(t,x;u(\cd))=\dbR\[\F(T,t)\xi+\int_t^T
\F(t,\t)
f(\t)d\t\],\q u(\cd)\in\sD(J(t,x;\cd)).$$
We see that $\sD(J(t,x;\cd))$ is a dense subspace of $\sU^2[t,T]$, but it does not have a clear description. To make such a description clearer, we introduce some  additional assumptions. To this end, similar to \rf{L^1+}, we introduce the following additional conditions.

\ms

{\bf(H4)} Let $G\in\dbS^n_+$, $g\in\dbR^n$, $q(\cd)\in L^{2+}(0,T;\dbR^n)$, and $r(\cd)\in L^{2+}(0,T;\dbR^m)$,
\bel{Q>0}\begin{pmatrix}Q(s)&S(s)^\top\\ S(s)&R(s)\end{pmatrix}\ges0,\qq s\in[0,T].\ee
Moreover,
\bel{R}\begin{pmatrix}q(s)\\ r(s)\end{pmatrix}\in\sR\(\begin{pmatrix}Q(s)&S(s)^\top\\ S(s)&R(s)\end{pmatrix}\),\qq g\in\sR(G).\ee

\ms

Under (H4), we have
$$g=GG^\dag g,$$
with $G^\dag$ being the pseudo-inverse of $G$. Then for any $u(\cd)\in\sU^2[t,T]$,
$$\ba{ll}
\ns\ds\xi=\lan GX(T),X(T)\ran+2\lan g,X(T)\ran=|G^{1\over2}X(T)|^2+2\lan G^{1\over2}G^\dag g,G^{1\over2}X(T)\ran\\
\ns\ds\q=\lan G[X(T)+G^\dag g],X(T)+G^\dag g\ran-\lan G^\dag g,g\ran\ges-\lan G^\dag g,g\ran\equiv\wt\xi.\ea$$
Similarly,
$$\ba{ll}
\ns\ds f(\cd)=\lan\begin{pmatrix}Q(s)&S(s)^\top\\ S(s)&R(s)\end{pmatrix}\begin{pmatrix}X(s)\\ u(s)\end{pmatrix},\begin{pmatrix}X(s)\\ u(s)\end{pmatrix}\ran+2\lan\begin{pmatrix}q(s)\\ r(s)\end{pmatrix},\begin{pmatrix}X(s)\\ u(s)\end{pmatrix}\ran\\
\ns\ds=\lan\begin{pmatrix}Q(s)&S(s)^\top\\ S(s)&R(s)\end{pmatrix}\[\begin{pmatrix}X(s)\\ u(s)\end{pmatrix}+\begin{pmatrix}Q(s)&S(s)^\top\\ S(s)&R(s)\end{pmatrix}^\dag\begin{pmatrix}q(s)\\ r(s)\end{pmatrix}\],\begin{pmatrix}X(s)\\ u(s)\end{pmatrix}+\begin{pmatrix}Q(s)&S(s)^\top\\ S(s)&R(s)\end{pmatrix}^\dag\begin{pmatrix}q(s)\\ r(s)\end{pmatrix}\ran\\
\ns\ds\q-\lan\begin{pmatrix}Q(s)&S(s)^\top\\ S(s)&R(s)\end{pmatrix}^\dag\begin{pmatrix}q(s)\\ r(s)\end{pmatrix},\begin{pmatrix}q(s)\\ r(s)\end{pmatrix}\ran\ges-\lan\begin{pmatrix}Q(s)&S(s)^\top\\ S(s)&R(s)\end{pmatrix}^\dag\begin{pmatrix}q(s)\\ r(s)\end{pmatrix},\begin{pmatrix}q(s)\\ r(s)\end{pmatrix}\ran\equiv\wt f(s),\q s\in[t,T].\ea$$
This implies $f(\cd)$ and $\xi$ bounded below by $\wt f(\cd)\in\h L^1_\dbF(\O;L^1(t,T;\dbR))$ and by $\wt\xi\in\h L^1_{\cF_T}(\O;\dbR)$. Therefore,
$$Y(s)=\dbE_s^t\[\F(T,s)\xi+\int_s^T\F(\t,s)
f(\t)d\t\],\qq s\in[t,T].$$
By comparison, if we set
$$\wt Y(s)=\dbE_s^t\[\F(T,s)\wt\xi+\int_s^T\F(\t,s)
\wt f(\t)d\t\],\qq s\in[t,T],$$
then
$$J(t,x;u(\cd))=Y(t)\ges\wt Y(t),\qq u(\cd)\in\sD(J(t,x;\cd)).$$
Since $\sU^p[t,T]$ ($p>2$) is dense in $\sU^2[t,T]$, we have
$$J(t,x;u(\cd))\ges\wt Y(t),\qq\forall u(\cd)\in\sD(J(t,x;\cd))).$$
Note that $Y(t)=\i$ is possible. Hence, under (H1), (H3), and (H4), we have
\bel{D*}\ba{ll}
\ns\ds\Big\{u(\cd)\1n\in\1n\sU^{2+}[t,T]\bigm|\hb{ \rf{f**}  holds}\Big\}\subseteq\Big\{u(\cd)\in\sU^2[t,T]\bigm|J(t,x;u(\cd))<\i\Big\}
\equiv\sD(J(t,x;\cd)).\ea\ee
and therefore, (note (H4))
\bel{2.30}\ba{ll}
\ns\ds\sU^2[t,T]=\Big\{u(\cd)\in\sU^2[t,T]\bigm|J(t,x;u(\cd))\les\i\Big\}\\
\ns\ds\qq\qq=\sD(J(t,x;u(\cd))
\bigcup\Big\{u(\cd)\in\sU^2[t,T]\bigm|J(t,x;u(\cd))=\i\Big\}.\ea\ee
Now, for $(t,x)\in[0,T]\times\dbR^n$ and $u(\cd)\in\sU^{2+}[t,T]$, we may let
$$X(s)=(\G^0x)(s)+\G^1(s)+\G_su(\cd)\equiv\begin{pmatrix}1&1&\G_s\end{pmatrix}
\begin{pmatrix}(\G^0x)(s)\\ \G^1(s)\\ \G_su(\cd)\end{pmatrix},\qq s\in[t,T].$$
Note that
$$\left\{\2n\ba{ll}
\ns\ds(\G^0x)(\cd)+\G^1(\cd)\in L^p_\dbF(\O;C([t,T];\dbR^n)),\qq p\ges1,\\
\ns\ds\G_\cd u(\cd)\in L^p_\dbF(\O;C([t,T];\dbR^n)),\qq u(\cd)\in\sU^p[t,T],\ea\right.$$
with estimate \rf{|X|}. Then one has
\bel{J=F}\ba{ll}
\ns\ds J(t,x;u(\cd))\\
\ns\ds=\dbE\[\F(T,t)
\lan G(1~1~\G_T)
\begin{pmatrix}(\G^0x)(T)\\ \G^1(T)\\
        u(\cd)\end{pmatrix},(1~1~\G_T)\begin{pmatrix}(\G^0x)(T)\\ \G^1(T)\\ u(\cd)\end{pmatrix}\ran\1n+\1n2\lan g,(1~1~\G_T)\begin{pmatrix}(\G^0x)(T)\\ \G^1(T)\\
u(\cd)\end{pmatrix}\ran\\
\ns\ds\qq+\int_t^T\F(\t,t)\(\lan\begin{pmatrix}Q(\t)&S(\t)^\top\\ S(\t)&R(\t)\end{pmatrix}
\begin{pmatrix}1&1&\G_\t\\
0&0&I\end{pmatrix}
\begin{pmatrix}(\G^0x)(\t)\\ \G^1(\t)\\
u(\cd)\end{pmatrix},\begin{pmatrix}1&1&\G_\t\\
0&0&I\end{pmatrix}
\begin{pmatrix}(\G^0x)(\t)\\ \G^1(\t)\\
u(\cd)\end{pmatrix}\ran\\
\ns\ds\qq\qq\qq+2\lan\begin{pmatrix}q(\t)\\ r(\t)\end{pmatrix},\begin{pmatrix}1&1&\G_\t\\
0&0&I\end{pmatrix}
\begin{pmatrix}(\G^0x)(\t)\\ \G^1(\t)\\
u(\cd)\end{pmatrix}\ran\)d\t\]\\
\ns\ds=\dbE\[\F(T,t)
\lan\begin{pmatrix}G&G&G\G_T\\ G&G& G\G_T\\ \G_T^*G&\G_T^*G&\G_T^*G\G_T\end{pmatrix}
\begin{pmatrix}(\G^0x)(T)\\ \G^1(T)\\
u(\cd)\end{pmatrix},
\begin{pmatrix}(\G^0x)(T)\\ \G^1(T)\\
u(\cd)\end{pmatrix}\ran+2\lan\begin{pmatrix}g\\ g\\ \G^*_Tg\end{pmatrix},
\begin{pmatrix}(\G^0x)(T)\\ \G^1(T)\\
u(\cd)\end{pmatrix}\ran\\
\ns\ds\qq+\int_t^T\F(\t,t)\(\lan\begin{pmatrix}Q(\t)&Q(\t)&Q(\t)\G_\t+S(\t)^\top\\
Q(\t)&Q(\t)&Q(\t)\G_\t+S(\t)^\top\\
\G_\t^*Q(\t)+S(\t)&\G_\t^*Q(\t)+S(\t)&\G_\t^*Q(\t)\G_\t+S(\t)^\top\G_\t+\G_\t^*S(\t) +R(\t)\end{pmatrix}\\
\ns\ds\qq\qq\begin{pmatrix}(\G^0x)(\t)\\ \G^1(\t)\\
u(\cd)\end{pmatrix},
\begin{pmatrix}(\G^0x)(\t)\\ \G^1(\t)\\
u(\cd)\end{pmatrix}\ran+2\lan\begin{pmatrix}q(\t)\\ q(\t)\\ \G_\t^*q(\t)+r(\t)\end{pmatrix},
\begin{pmatrix}(\G^0x)(\t)\\ \G^1(\t)\\
u(\cd)\end{pmatrix}\ran\)d\t\]\\
\ns\ds\equiv\lan\F_2u(\cd),u(\cd)\ran+2\lan\F_1,u(\cd)\ran+\F_0,\ea\ee
where
$$\ba{ll}
\ns\ds\lan\F_2u(\cd),u(\cd)\ran=\dbE\[\F(T,t)\lan G
\G_Tu(\cd),\G_Tu(\cd)\ran\\
\ns\ds\qq\qq\qq\qq+\int_t^T\F(\t,t)\(\lan Q(\t)\G_\t u(\cd),\G_\t u(\cd)\ran+2\lan S(\t)\G_\t u(\cd),u(\t)\ran+\lan R(\t)u(\t),u(\t)\ran\)d\t\],\\
\ns\ds\lan\F_1,u(\cd)\ran=\dbE\[\lan g,\G_Tu(\cd)\ran+\int_t^T\F(\t,t)\(\lan(\G^0x)(\t)+\G^1(\t)+r(\t),u(\t)\ran+\lan q(\t),\G_\t u(\cd)\ran\)d\t\],\\
\ns\ds\F_0=\dbE\[\F(T,t)
\lan\begin{pmatrix}G&G\\ G&G\end{pmatrix}
\begin{pmatrix}(\G^0x)(T)\\ \G^1(T)\end{pmatrix},
\begin{pmatrix}(\G^0x)(T)\\ \G^1(T)\end{pmatrix}\ran+2\lan\begin{pmatrix}g\\ g\end{pmatrix},
\begin{pmatrix}(\G^0x)(T)\\ \G^1(T)\end{pmatrix}\ran\\
\ns\ds\qq+\int_t^T\F(\t,t)\(\lan\begin{pmatrix}Q(\t)&Q(\t)\\
Q(\t)&Q(\t)\end{pmatrix}\begin{pmatrix}(\G^0x)(\t)\\ \G^1(\t)\end{pmatrix},
\begin{pmatrix}(\G^0x)(\t)\\ \G^1(\t)\end{pmatrix}\ran+2\lan\begin{pmatrix}q(\t)\\ q(\t)\end{pmatrix},
\begin{pmatrix}(\G^0x)(\t)\\ \G^1(\t)\end{pmatrix}\ran\)d\t\].\ea$$
From this, we have the following result.

\bl{Lemma 2.4} \sl Let {\rm(H1)}, {\rm(H3)}, and
{\rm(H4)} hold. Then for any $t\in[0,T]$, $(x,u(\cd))\mapsto J(t,x;u(\cd))$ is a quadratic form and convex on $\sD(J(t,x;\cd))$.

\el

\it Proof. \rm Clearly, for any $(t,x)\in[0,T]\times\dbR^n$ and $u(\cd)\in\sU^p[t,T]$, with $p>2$, the map
\bel{X=X+X+X}(t,x,u(\cd))\mapsto X(\cd\,;t,x,u(\cd))\equiv X(\cd)\equiv (\G^0x)(\cd)+\G^1(\cd)+\G_\cd u(\cd),\ee
is affine and
$$J(t,x;u(\cd))=\lan\F_2u(\cd),u(\cd)\ran+2\lan\F_1,u(\cd)\ran+\F_0$$
is quadratic and convex on $\sU^2[t,T]$. Since $\sU^p[t,T]$ is dense in $\sD(J(t,x;\cd))$, our conclusion follows. \endpf

\ms

Now, we introduce the following notions.

\bde{solvability} Consider Problem (LQR), under (H1), (H3), (H4),

\ms

(i) Problem (LQR) is open-loop solvable at $(t,x)\in[0,T]\times\dbR^n$, if there exists an {\it open-loop optimal control} $\bar u(\cd)\in\sD(J(t,x;\cd))$ such that \rf{bar u} is satisfied.

\ms

(ii) A pair of maps $\Psi:[t,T]\to\dbR^m$ (bounded) and $v(\cd)\in\sU^2[t,T]$ is called an {\it affine strategy} if the following closed-loop system is uniquely solvable for any $x\in\dbR^n$:
\bel{closed-loop system}\left\{\2n\ba{ll}
\ns\ds dX(s)=\big\{[A(s)+B(s)\Psi(s)]X(s)+B(s)v(s)+b(s)\big\}ds\\
\ns\ds\qq\qq+\big\{[C(s)+D(s)\Psi(s)]X(s)+D(s)v(s)+\si(s)\big\}dW(s),\q s\in[t,T],\\
\ns\ds X(t)=x.\ea\right.\ee
Let $\BBPsi[t,T]\times\sU^2[t,T]$ be the set of all affine feedback strategy on $[t,T]$. In this case,
\bel{u()}u(\cd)=\Psi(\cd)X(\cd)+v(\cd)\in\sD(J(t,x;\cd))\ee
is called the {\it outcome} of the strategy $(\Psi(\cd),v(\cd))$. The corresponding cost is denoted by
$$J(t,x;\Psi(\cd),v(\cd))=J(t,x;\Psi(\cd)X(\cd)+v(\cd)),$$
with $X(\cd)$ being the corresponding state process of the closed-loop system.

\ms

(iii) Problem (LQR) is closed-loop solvable at $t\in[0,T]$ if there exists a pair $(\bar\Psi(\cd),\bar v(\cd))\in\BBPsi[t,T]\times\sU^2[t,T]$ so that
\bel{J=inf*}J(t,x;\bar\Psi(\cd),\bar v(\cd))=\inf_{(\Psi(\cd),v(\cd))\in\BBPsi[t,T]\times\sU^2[t,T]}
J(t,x;\Psi(\cd),v(\cd)),\qq\forall x\in\dbR^n.\ee

\ede

\ms

Note that on the right-hand side of \rf{J=inf*} is either $\i$, for which $\Psi(\cd)X(\cd)+v)(cd)\notin\sD(J(t,x;\cd))$ and the strategy $(\Psi(\cd),v(\cd))$ will not be taken, or finite, for which $\Psi\Psi(\cd)X(\cd)+v(\cd)\in\cD(J(t,x;cd))$ and the strategy $(\Psi(\cd),v(\cd))$ might be taken. Also under (H1), (H3), and (H4), \rf{f^-xi^-} is always true. Thus, due to \rf{2.30}, in the above definition, $\sD(J(t,x;\cd))$ and
$\sU^2[t,T]$ are exchangeable. We point out that since for different initial state $x\in\dbR^n$, the state process $X(\cd)$ might be different. Hence, the outcome should be different. Consequently, \rf{J=inf*} actually gives a family of inequalities. Hence, if Problem (LQR) is open-loop solvable at $(t,x)\in[0,T]\times\dbR^n$, for each $x\in\dbR^n$ with the outcome $u(\cd)\in
\sD(J(t,x;\cd))$ of the same strategy $(\Psi(\cd),v(\cd))\in\BBPsi[t,T]\times\sU^2[t,T]$, then Problem (LQR) is closed-loop solvable at $t$. Finally, to conclude this section, we present the following:

\bl{Lemma 2.5} \sl Let {\rm(H1), (H3)} and {\rm(H4)} hold. Then $(\bar\Psi(\cd),\bar v(\cd))$ is an optimal closed-loop strategy on $[t,T]$ if and only if
\bel{J(Psi)}J(t,x;\bar\Psi(\cd)\bar X(\cd)+\bar v(\cd))\les J(t,x;u(\cd)),\qq\forall u(\cd)\in\sU^2[t,T].\ee

\el

\it Proof. \rm Suppose \rf{J=inf*} holds. Then for any $u(\cd)\in\sU^2[t,T]$, by taking $(\Psi(\cd),v(\cd))=(0,u(\cd))$, we have \rf{J(Psi)}. Conversely, let \rf{J(Psi)} holds. Then, for any affine strategy $(\Psi(\cd),v(\cd))\in\BBPsi[t,T]\times\sU^2[t,T]$, let
$$u(\cd)=\Psi(\cd)X(\cd)+v(\cd),$$
with $X(\cd)$ being the solution of \rf{closed-loop system}. Then \rf{J(Psi)} implies \rf{J=inf*}. \endpf

\ms

The above result means if $(\bar\Psi(\cd),\bar v(\cd))$ is a closed-loop optimal strategy, then its outcome $\bar\Psi(\cd)\bar X(\cd)+\bar v(\cd)$ is an open-loop optimal control. Hence, Problem (LQR) is closed-loop solvable on $[t,T]$, then it is open-loop solvable at any $(t,x)\in[0,T]\times\dbR^n$. We know that, in general, the converse is not true (see \cite{Sun-Yong 2020}, Example 2.1.6).

\ms

\section{A Relation to the Classical Case.}

In this section, we are transforming Problem (LQR) to a classical one, and reveal a relation of our Problem (LQR) with the classical Problem (LQ). First of all, let
\bel{L}\L(s,t)=e^{\int_t^s[-{1\over2}F(\th)^2+E(\th)]d\th},\qq0\les s\les t\les T.\ee
Thus, the following type group property holds:
$$\L(s,t)=\L(s,\t)\L(\t,t),\qq0\les s,\t,t\les T.$$
When $F(\cd)$ and $E(\cd)$ are bounded, $\L(\cd\,,\cd)$ is bounded as well. Then
\bel{Y(s)}\ba{ll}
\ns\ds Y(s)=\dbE_s^t\[\F(T,s)\(\xi+\int_s^T\F(\t,T)
e^{\int_T^\t E(\th)d\th}f(\t)d\t\)\]\\
\ns\ds\qq=\dbE^t_s\[e^{\int_s^TF(\th)dW(\th)+\int_s^T[-{1\over2}F(\th)^2+E(\th)]
d\th}\(\xi+\int_s^Te^{\int_T^\t F(\th)dW(\th)+\int_T^\t[-{1\over2}F(\th)^2+E(\th)]d\th}f(\t)d\t\)\]\\
\ns\ds\qq=\dbE^t_s\[e^{\int_s^TF(\th)dW(\th)}\L(T,s)\(\xi+\int_s^Te^{\int_T^\t F(\th)dW(\th)}\L(\t,T)f(\t)d\t\)\]\\
\ns\ds\qq=\L(T,s)\dbE_s^t\(e^{\int_s^TF(\th)dW(\th)}\xi+\int_s^Te^{\int_s^\t F(\th)dW(\th)}\L(\t,T)f(\t)d\t\).\ea\ee
Thus,
\bel{J*}\ba{ll}
\ns\ds J(t,x;u(\cd))=Y(t)=\L(T,t)\dbE\[e^{\int_t^TF(\th)dW(\th)}\(\lan GX(T),X(T)\ran+2\lan g,X(T)\ran\)\\
\ns\ds\qq\qq\qq\qq\qq+\int_t^T\L(\t,T)e^{\int_t^\t F(\th)dW(\th)}\(\lan Q(\t)X(\t),X(\t)\ran+2\lan S(\t)X(\t),u(\t)\ran\\
\ns\ds\qq\qq\qq\qq\qq+\lan R(\t)u(\t),u(\t)\ran+2\lan q(\t),X(\t)\ran+2\lan r(\t),u(\t)\ran\)d\t\]\\
\ns\ds=\L(T,t)\dbE\[\lan G\h X(T),\h X(T)\ran+2\lan\h g,\h X(T)\ran+\int_t^T\L(\t,T)\(\lan Q(\t)\h X(\t),\h X(\t)\ran\\
\ns\ds\qq+2\lan S(\t)\h X(\t),\h u(\t)\ran+\lan R(\t)\h u(\t),\h u(\t)\ran+2\lan \h q(\t),\h X(\t)\ran+2\lan\h r(\t),\h u(\t)\ran\)d\t\]\ea\ee
where
\bel{trans}\left\{\2n\ba{ll}
\ds\h g=e^{{1\over2}\int_t^TF(\th)dW(\th)}g,\q\h q(\t)=e^{\int_t^\t{1\over2}F(\th)dW(\th)}q(\t),\q\h r(\t)=e^{\int_t^\t{1\over2}F(\th)dW(\th)}r(\t),\\
\ns\ds\h X(\t)=e^{{1\over2}\int_t^\t F(\th)dW(\th)}X(\t),\q\h u(\t)=e^{{1\over2}\int_t^\t F(\th)dW(\th)}u(\t),\ea\right.\ee
Applying the It\^o's formula to the map $\t\mapsto e^{{1\over2}\int_t^\t F(\th)dW(\th)}X(\t)$, one has
$$\ba{ll}
\ns\ds d\h X(\t)=\[e^{{1\over2}\int_t^\t F(\th)dW(\th)}{1\over2}F(\t)dW(\t)\]X(\t)+e^{{1\over2}\int_t^\t F(\th)
dW(\th)}\([A(\t)X(\t)+B(\t)u(\t)+b(\t)]d\t\\
\ns\ds\qq\qq+[C(\t)X(\t)+D(\t)u(\t)+\si(\t)]dW(\t)\)
+{1\over2}F(\t)e^{{1\over2}\int_t^\t F(\th)
dW(\th)}[C(\t)X(\t)+D(\t)u(\t)+\si(\t)]d\t\\
\ns\ds\qq\q={1\over2}\h X(\t)F(\t)dW(\t)+[A(\t)\h X(\t)+B(\t)\h u(\t)+e^{{1\over2}\int_t^\t F(\th)dW(\th)}b(\t)]d\t\\
\ns\ds\qq\qq+[C(\t)\h X(\t)+D(\t)\h u(\t)+e^{{1\over2}\int_t^\t F(\th)dW(\th)}\si(\t)]dW(\t)\\
\ns\ds\qq\qq+{1\over2}F(\t)[C(\t)\h X(\t)+D(\t)\h u(\t)+e^{{1\over2}\int_t^\t F(\th)dW(\th)}\si(\t)]d\t\\
\ns\ds\qq\q=\([A(\t)+{1\over2}F(\t)C(\t)]\h X(\t)+[B(\t)+{1\over2}F(\t)D(\t)]\h u(\t)+e^{{1\over2}\int_t^\t F(\th)dW(\th)}[b(\t)+{1\over2}F(\t)\si(\t)]\)d\t\\
\ns\ds\qq\qq+\([C(\t)+{1\over2}F(\t)I]\h X(\t)+D(\t)\h u(\t)+e^{{1\over2}\int_t^\t F(\th)dW(\th)}\si(\t)\)dW(\t)\ea$$
Then, we have a classical LQ problem with the state equation:
\bel{hX}\left\{\2n\ba{ll}
\ds d\h X(\t)=\(\h A(\t)\h X(\t)+\h B(\t)\h u(\t)+\h b(\t)\)d\t\\
\ns\ds\qq\qq+\(\h C(\t)\h X(\t)+\h D(\t)\h u(\t)+\h\si(\t)\)dW(\t),\q\t\in[t,T],\\
\ns\ns\h X(t)=x,\ea\right.\ee
and with cost functional:
\bel{wt J}\ba{ll}
\ns\ds\h J(t,x;\h u(\cd))=\dbE\[\(\lan G\h X(T),\h X(T)\ran+2\lan\h g,\h X(T)\ran\)\\
\ns\ds\qq\qq\qq\qq+\int_t^T\L(\t,T)\(\lan  Q(\t)\h X(\t),\h X(\t)\ran+2\lan S(\t)\h X(\t),\h u(\t)\ran\\
\ns\ds\qq\qq\qq\qq\qq+\lan R(\t)\h u(\t),\h u(\t)\ran+2\lan\h q(\t),\h X(\t)\ran+2\lan\h r(\t),\h u(\t)\ran\)d\t\].\ea\ee
In the above, we have used
\bel{h A}\left\{\2n\ba{ll}
\ds\h A(\t)=A(\t)+{1\over2}F(\t)C(\t),\qq\h B(\t)=B(\t)+{1\over2}F(\t)D(\t),\\
\ns\ds\h C(\t)=C(\t)+{1\over2}F(\t)I,\qq\h D(\t)=D(\t),\\
\ns\ds\h b(\t)=e^{{1\over2}\int_t^\t F(\th)dW(\th)}\(b(\t)+{1\over2}F(\t)\si(\t)\),\qq\h\si(\t)=e^{{1\over2}\int_t^\t F(\th)dW(\th)}\si(\t).\ea\right.\ee
Also,
\bel{J}J(t,x;u(\cd))=Y(t)=\L(T,t)\h J(t,x;\h u(\cd)).\ee
Hence, under transformation \rf{trans}, we obtain the following classical LQ problem.

\ms

\bf Problem (LQ). \rm For given $(t,x)\in[0,T]\times\dbR^n$, find $\cl{\h u}(\cd)\in\sU^2[t,T]$ such that
\bel{J=inf}\h J(t,x;\cl{\h u}(\cd))=\inf_{\h u(\cd)\in\sU^2[t,T]}\h J(t,x;\h u(\cd)).\ee

\ms

Problem (LQR) seems to be transformed to the above classical Problem (LQ). But, we might overlook one point. Function $u(\cd)\mapsto J(t,x;u(\cd))$ is defined on $\sD(J(t,x;\cd))$ which is dense in $\sU^2[t,T]$, whereas $u(\cd)\mapsto\h J(t,x;u(\cd))$ is defined on $\sU^2[t,T]$. Thus, in some sense, $\h J(t,x;\cd)$ is an extension of $J(t,x;\cd)$. Now, we have two Problems (LQ), one is with functional defined on $\sD(J(t,x;\cd))$ (transformed from Problem (LQR)) and the other is the classical one with the cost defined on
$\sU^2[t,T]$. It is known that, under some mild conditions, the latter has a minimum $\cl{\h u}(\cd)\in\sU^2[t,T]$. But $\cl{\h u}(\cd)$ might not be in $\sD(J(t,x;\cd))$ and therefore it could not be a minimum of $J(t,x;u(\cd))$.
To make sure if $\cl{\h u}(\cd)$ is in $\sD(J(t,x;\cd))$ or not is a hard job. On the other hand, $\h b(\cd)$ and $\h\si(\cd)$ defined in \rf{h A} are random, which will lead to some additional technicality. Thus, we prefer to approach the problem directly.

\section{FBSDEs and Open-Loop Solvability}

In this section, we are going to look at the open-loop solvability of Problem (LQR). Let us first present a general result, which is similar to that in \cite{Mou-Yong 2006}.

\bl{L4.1} \sl Let $\BH$ be a Hilbert space and $\BF:\sD(F)\subseteq\BH\to\dbR$ be a quadratic function, i.e.,
\bel{F}\BF(u)=\lan\F_2u,u\ran+2\lan\F_1,u\ran+\F_0,\qq u\in\sD(\BF)\equiv\sD(\F_2),\ee
with $\F_2:\sD(\F_2)\subseteq\BH\to\BH$ being self-adjoint, $\F_1\in\BH$ and $\F_0\in\dbR$. Let $u\in\sD(\F_2)$ attend a minimum of $\BF$. Then
\bel{F_2>0}\F_2\ges0,\ee
and
\bel{F_2u+F_1=0}D_u\BF(u)=2\big(\F_2u+\F_1\big)=0,\ee
admits a solution $\bar u\in\sD(\F_2)$. Conversely. if \rf{F_2>0} holds and $\bar u\in\sD(\F_2)$ solves equation \rf{F_2u+F_1=0}, then $\bar u$ is a minimum of $\BF$ over $\sD(\F_2)\subseteq\BH$. In this case.
\bel{F*}\BF(u)=\lan\F_2(u-\bar u),u-\bar u\ran+\F_0-\lan\F_2^\dag\F_1,\F_1\ran,\ee
where $\F_2^\dag$ is the pseudo-inverse of $\F_2\ges0$. In the case that $\F_2$ is uniformly positive definite, the minimum uniquely exists.

\el

\it Proof. \rm Let $\bar u\in\sD(\F_2)\subseteq\BH$ such that
$$\BF(\bar u)=\min_{u\in\sD(\BF)\subseteq\BH}\BF(u)=\min_{u\in\sD(\F_2)\subseteq\BH}
\[\lan\F_2u,u\ran+2\lan\F_1,u\ran+\F_0\]>-\i.$$
Then we claim \rf{F_2>0} holds. In fact, if it were not true. Then, there exists a $u_0\in\sD(\F_2)$ such that
$$\lan\F_2u_0,u_0\ran<0.$$
Then,
$$\BF(\l u_0)=\l^2\lan\F_2u_0,u_0\ran+2\l\lan\F_1,u_0\ran+\F_0\to-\i,$$
which is a contradiction. Hence, \rf{F_2>0} is true.

\ms

Next, suppose $\BF(\cd)$ attends a minimum at $\bar u\in\sD(\F_2)\subseteq\BH$. Then, by Fermat' Theorem, we see that $\bar u$ is a solution of equation \rf{F_2u+F_1=0}. Conversely, if \rf{F_2>0} holds and $\bar u\in\sD(\F_2)\subseteq\BH$ solves equation \rf{F_2u+F_1=0}, then (note $\F_1=-\F_2\bar u$)
\bel{F(u)}\ba{ll}
\ns\ds\BF(u)=\lan\F_2u,u\ran+2\lan\F_1,u\ran+\F_0=|\F_2^{1\over2}u|^2-2\lan\F_2\bar u,u\ran+\F_0\\
\ns\ds\qq=|\F_2^{1\over2}u|^2-2\lan\F_2^{1\over2}u,\F_2^{1\over2}\bar u\ran+\F_0=|\F_2^{1\over2}(u-\bar u)|^2+\F_0-|\F_2^{1\over2}\bar u|^2\\
\ns\ds\qq=\lan\F_2(u-\bar u),u-\bar u\ran+\F_0-\lan\F_2\bar u,\bar u\ran
=\lan\F_2(u-\bar u),u-\bar u\ran+\F_0-\lan\F_2^\dag\F_1,\F_1\ran.\ea\ee
since $\F_2\F_2^\dag\F_2=\F_2$. Now, if $\F_2$ is uniformly positive definite.
Then for any minimizing sequence $u_k$,
$$\ba{ll}
\ns\ds\BF(u_k)=\lan\F_2u_k,u_k\ran+2\lan\F_1,u_k\ran+\F_0=|\F_2^{1\over2}u_k|^2
+\lan\F_2^{1\over2}u_k,\F_2^{-{1\over2}}\F_1\ran+\F_0\\
\ns\ds\qq\q=|\F_2^{1\over2}(u_k-\F_2^{-1}\F_1)|^2+\F_0-\lan\F_2^{-1}\F_1,\F_1\ran.\ea$$
Hence, $u_k$ is bounded. Then, we may assume $u_k$ goes to some $\bar u$ weekly. Then by the lower semi-continuity of $\BF(\cd)$, we have
$$\BF(\bar u)\les\liminf_{k\to\i}\BF(u_k).$$
Thus, $\bar u$ is a minimum of $\BF$. Thus, by \rf{F(u)},
$$\min_{u\in\BH}\BF(u)=\F_0-\lan\F_2^{-1}\F_1,\F_1\ran=\BF(u')=\lan\F_2(u'-\bar u),u'-\bar u\ran+\F_0-\lan\F_2^{-1}\F_1,\F_1\ran,$$
leading to
$$\d|u'-\bar u|^2\les\lan\F_2(u'-\bar u),u'-\bar u\ran=0.$$
This gives the uniqueness.  \endpf

\ms

Next, from \rf{F} and \rf{J=F}, we see that for a given state-control pair $(X(\cd),u(\cd))$,
$$\lan D_{u(\cd)}J(t,x;u(\cd)),u_1(\cd)\ran=2\lan\F_2u+\F_1,u_1\ran.$$
In the above, $\F_2$ and $\F_1$ are a little too abstract.
We have the following result, which gives the Fr\'echet derivative of $u(\cd)\mapsto J(t,x;u(\cd))$.

\bt{Gateaux} \sl Let {\rm(H1)}, {\rm(H3)}, and {\rm(H4)} hold. Let $(t,x)\in[0,T]\times\dbR^n$ and $u(\cd)\in\sD(J(t,x;\cd))$. Let $(X(\cd),u(\cd))$ be a state-control pair. Then, for any $u_1(\cd)\in\sD(J(t,x;\cd))$,
the G\^ateaux derivative of $u(\cd)\mapsto J(t,x;u(\cd))$ is given by
\bel{D_uJ}\ba{ll}
\ns\ds\lan D_{u(\cd)}J(t,x;u(\cd)),u_1(\cd)\ran=2\dbE\int_t^T
\F(s,t)\(\lan[B(s)+F(s)D(s)]^\top\wt Y(s)+D(s)^\top\wt Z(s)\\
\ns\ds\qq\qq\qq\qq\qq\qq\qq\qq\qq+S(s)X(s)+R(s)u(s)+r(s),u_1(s)\ran\)ds,\ea\ee
with $(\wt Y(\cd),\wt Z(\cd))$ being the adapted solution to following $n$-dimensional BSDE:
\bel{wt Y}\left\{\2n\ba{ll}
\ds d\wt Y(s)=-\[\(A(s)+E(s)I+F(s)C(s)\)^\top\wt Y(s)+\(C(s)+F(s)I\)^\top\wt Z(s)\\
\ns\ds\qq\qq\qq\qq\qq\qq+Q(s)X(s)+S(s)^\top u(s)+q(s)\]ds+\wt Z(s)dW(s),\q s\in[t,T],\\
\ns\ds\wt Y(T)=GX(T)+g.\ea\right.\ee
\et

\it Proof. \rm For given $(t,x)\in[0,T]\times\dbR^n$, suppose $(X(\cd),u(\cd))$ is a state-control pair. Pick any $u_1(\cd)\in\sU^p[t,T]$, and $\e>0$ small. Let $X^\e(\cd)$ be the state process corresponding to $(t,x)$ and
$$u^\e(\cd)=u(\cd)+\e u_1(\cd)\in\sU^2[t,T].$$
Then
$${X^\e(\cd)-X(\cd)\over\e}\to X_1(\cd),\qq\hb{in }L^1_\dbF(\O;C([t,T];\dbR^n)),$$
with $X_1(\cd)$ being the solution of the following variational system:
\bel{variation}\left\{\2n\ba{ll}
\ds dX_1(s)=[A(s)X_1(s)+B(s)u_1(s)]ds+[C(s)X_1(s)+D(s)u_1(s)]dW(s),\q s\in[t,T],\\
\ns\ds X_1(t)=0.\ea\right.\ee
Since $u_1(\cd)\in\sU^p[t,T]$ (with $p>2$), $X_1(\cd)\in L^p_\dbF(\O;
C([t,T];\dbR^n))$, and
\bel{|X_1|}\dbE\[\sup_{s\in[t,T]}|X_1(s)|^p\]\les K\dbE\(\int_t^T|u(s)|^2ds\)^{p\over2}.\ee
Now, we let
$$Y(s)=\dbE^t_s\[\F(T,s)\xi(t,x; X(T))+\int_s^T\F(\t,s)f(t,x;X(\t),u(\t))d\t\],$$
with $\xi$ and $f(\cd)$ given by \rf{f**}, and
$$Y^\e(s)=\dbE^t_s\[\F(t,s)\xi(t,x; X^\e(T))+\int_s^T\F(\t,s)f(t,x;X^\e(\t),u^\e(\t))d\t\].$$
Hence, we may let
$${Y^\e(\cd)-Y(\cd)\over\e}\to\h Y_1(\cd),\q\hb{in }L^1_\dbF(\O;C([t,T];\dbR^n)),$$
with $\h Y_1(\cd)$ being the solution of the following variational system:
$$\ba{ll}
\ns\ds\h Y_1(s)=\dbE_s^t\[\F(T,s)\(2\lan GX(T)+g,X_1(T)\ran\)+\int_s^T\F(\t,s)\(2\lan Q(\t)X(\t)+S(\t)^\top u(\t)+q(\t),X_1(\t)\ran\\
\ns\ds\qq\qq\qq\qq+2\lan S(\t)X(\t)+R(\t)u(\t)+r(\t),u_1(\t)\ran\)
d\t\]\\
\ns\ds\qq\equiv\dbE_s^t\[\F(T,s)\xi_1+\int_s^T
\F(\t,s)f_1(\t)d\t\].\ea$$
Since $u_1(\cd)\in\sU^p[t,T]$, with $p>2$, we have \rf{|X_1|}.
Consequently, in the above, $(\xi_1,f_1(\cd))\in L^{p'}_{\cF_T}(\O;\dbR)\times L^{p'}_\dbF(\O;L^{p'}(t,T;\dbR))$, for
some $p'>2$. We see that the above is equivalent to the following
1-dimensional classical BSDE:
$$\left\{\2n\ba{ll}
\ds d\h Y_1(s)=-[E(s)\h Y_1(s)+F(s)\h Z_1(s)+2\lan Q(s)X(s)+S(s)^\top u(s)+q(\t),X_1(s)\ran\\
\ns\ds\qq\qq\qq+2\lan S(s)X(s)+R(s)u(s)+r(s),u_1(s)\ran ds+\h Z_1(s)dW(s),\q s\in[t,T],\\
\ns\ds\h Y_1(T)=2\lan GX(T)+g,X_1(T)\ran,\ea\right.$$
for which $(\h Y_1(\cd),\h Z_1(\cd))$ is the adapted solution. Next, we let
$(\wt Y(\cd),\wt Z(\cd))$ be the adapted solution of the following $n$-dimensional BSDE:
$$\left\{\2n\ba{ll}
\ds d\wt Y(s)=-\wt\G(s)ds+\wt Z(s)dW(s),\qq s\in[t,T],\\
\ns\ds\wt Y(T)=GX(T)+g,\ea\right.$$
with $\wt\G(\cd)$ undetermined. Note that since the terminal state is
in $L^2_{\cF_T}(\O;\dbR^n)$. The above BSDE is well-posed (provided the
undetermined $\wt\G(\cd)$ is in $L^{p\over2}_\dbF(\O;L^{p\over2}(t,T;
\dbR))$ for some $p>2$). Then by It\^o's formula,
$$\ba{ll}
\ns\ds d\lan\wt Y(s),X_1(s)\ran=[\lan-\wt\G(s)ds+\wt Z(s)dW(s),X_1(s)\ran+\lan\wt Y(s),[A(s)X_1(s)+B(s)u_1(s)]ds\\
\ns\ds\qq\qq\qq\qq+[C(s)X_1(s)+D(s)u_1(s)]dW(s)\ran+\lan\wt Z(s),C(s)X_1(s)+D(s)u_1(s)\ran ds\\
\ns\ds\qq\qq\qq\q=[\lan-\wt\G(s)+A(s)^\top\wt Y(s)+C(s)^\top\wt Z(s),X_1(s)\ran
+\lan B(s)^\top\wt Y(s)+D(s)^\top\wt Z(s),u_1(s)\ran]ds\\
\ns\ds\qq\qq\qq\qq+[\lan\wt Z(s),X_1(s)\ran+\lan\wt Y(s),C(s)X_1(s)+D(s)u_1(s)\ran]dW(s).\ea$$
Now, let
$$\left\{\2n\ba{ll}
\ds Y_1(s)=\h Y_1(s)-2\lan\wt Y(s),X_1(s)\ran,\\
\ns\ds Z_1(s)=\h Z_1(s)-2[\lan\wt Z
(s),X_1(s)\ran+\lan\wt Y(s),C(s)X_1(s)+D(s)u_1(s)\ran].\ea\right.$$
Then
\bel{Y_1}Y_1(T)=\h Y_1(T)-2\lan\wt Y(T),X_1(T)\ran=2\lan GX(T)+g,X_1(T)\ran-2\lan GX(T)+g,X_1(T)\ran=0,\ee
and
$$\ba{ll}
\ns\ds dY_1(s)=d\h Y_1(s)-2d\lan\wt Y(s),X_1(s)\ran\\
\ns\ds\qq\q=-[E(s)\h Y_1(s)+F(s)\h Z_1(s)+2\lan Q(s)X(s)+S(s)^\top u(s)+q(s),X_1(s)\ran\\
\ns\ds\qq\qq\qq+2\lan S(s)X(s)+R(s)u(s)+r(s),u_1(s)\ran]ds+\h Z_1(s)dW(s)\\
\ns\ds\qq\qq-2[\lan-\wt\G(s)+A(s)^\top\wt Y(s)+C(s)^\top\wt Z(s),X_1(s)\ran
+\lan B(s)^\top\wt Y(s)+D(s)^\top\wt Z(s),u_1(s)\ran]ds\\
\ns\ds\qq\qq\qq\qq-2[\lan\wt Z(s),X_1(s)\ran+\lan\wt Y(s),C(s)X_1(s)+D(s)u_1(s)\ran]dW(s)\ea$$
$$\ba{ll}
\ns\ds\qq\q=-\[E(s)\(Y_1(s)+2\lan\wt Y(s),X_1(s)\ran\)\\
\ns\ds\qq\qq\qq+F(s)\(Z_1(s)+2[\lan\wt Z(s),X_1(s)\ran+\lan\wt Y(s), C(s)X_1(s)+D(s)u_1(s)\ran]\)\\
\ns\ds\qq\qq\qq+2\lan-\wt\G(s)+A(s)^\top\wt Y(s)+C(s)^\top\wt Z(s)+Q(s)X(s)+S(s)^\top u(s)+q(s),X_1(s)\ran\\
\ns\ds\qq\qq\qq+2\lan B(s)^\top\wt Y(s)+D(s)^\top\wt Z(s)+S(s) X(s)+R(s)u(s)+r(s),u_1(s)\ran\]ds\\
\ns\ds\qq\qq\qq\qq+\[\h Z_1(s)-2[\lan\wt Z(s),X_1(s)\ran+\lan\wt Y(s),C(s)X_1(s)+D(s)u_1(s)\ran\]dW(s)\\
\ns\ds\qq\q=-\[E(s)Y_1(s)+F(s)Z_1(s)+2\lan-\wt\G(s)+E(s)\wt Y(s)+F(s)\wt Z(s)+F(s)C(s)^\top\wt Y(s)\\
\ns\ds\qq\qq\qq+A(s)^\top\wt Y(s)+C(s)^\top\wt Z(s)+Q(s)X(s)+S(s)^\top u(s)+q(s),X_1(s)\ran\\
\ns\ds\qq\qq\qq+2\lan F(s)D(s)^\top\wt Y(s)+
B(s)^\top\wt Y(s)+D(s)^\top\wt Z(s)+S(s)X(s)+R(s)u(s)+r(s),u_1(s)\ran\]ds\\
\ns\ds\qq\qq\qq+Z_1(s)dW(s)\\
\ns\ds\qq\q=-\[E(s)Y_1(s)+F(s)Z_1(s)+2\lan-\wt\G(s)+[A(s)^\top+E(s)I+F(s)C(s)^\top]\wt Y(s)\\
\ns\ds\qq\qq\qq+[C(s)^\top+F(s)I]\wt Z(s)+Q(s)X(s)+S(s)^\top u(s)+q(s),X_1(s)\ran\\
\ns\ds\qq\qq\qq+2\lan B(s)^\top+F(s)D(s)^\top]\wt Y(s)+D(s)^\top\wt Z(s)+S(s)X(s)+R(s)u(s)+r(s),u_1(s)\ran\]ds\\
\ns\ds\qq\qq\qq+Z_1(s)dW(s).\ea$$
Thus, we take
$$\ba{ll}
\ns\ds\wt\G(s)=[A(s)+E(s)I+F(s)C(s)^\top]\wt Y(s)+[C(s)+F(s)I]^\top\wt Z(s)+Q(s)X(s)+S(s)^\top u(s)+q(s).\ea$$
Hence, we have the following BSDEs:
$$\left\{\2n\ba{ll}
\ds d\wt Y(s)=-\([A(s)+E(s)I+F(s)C(s)]^\top\wt Y(s)+[C(s)+F(s)I]^\top\wt Z(s)\\
\ns\ds\qq\qq\qq\qq\qq\qq\qq+Q(s)X(s)+S(s)^\top u(s)+q(s)\)ds+\wt Z(s)dW(s),\qq s\in[t,T],\\
\ns\ds\wt Y(T)=GX(T)+g,\ea\right.$$
and
$$\left\{\2n\ba{ll}
\ds dY_1(s)=-\(E(s)Y_1(s)+F(s)Z_1(s)+2\lan[B(s)+F(s)D(s)]^\top\wt Y^u(s)+D(s)^\top\wt Z^u(s)\\
\ns\ds\qq\qq\qq+S(s)X(s)+R(s)u(s)+r(s),u_1(s)\ran\)ds+Z_1(s)dW(s),\q s\in[t,T],\\
\ns\ds Y_1(T)=0.\ea\right.$$
Then noting \rf{Y_1} and $X_1(t)=0$, we have
$$\ba{ll}
\ns\ds{J(t,x;u^\e(\cd))-J(t,x;u(\cd))\over\e}={Y^\e(t)-Y(t)\over\e}\to\lan D_{u(\cd)}J(t,x;u(\cd)),u_1(\cd)\ran=\h Y_1(t)=Y_1(t)\\
\ns\ds=2\dbE\int_t^T
\F(s,t)\(\lan[B(s)+F(s)D(s)]^\top\wt Y(s)+D(s)^\top\wt Z(s)\\
\ns\ds\qq\qq\qq\qq\qq\qq\qq\qq\qq+S(s)X(s)+R(s)u(s)+r(s),u_1(s)\ran\)ds,\ea$$
Since $\sU^{2+}[t,T]$ is dense in $\sD((J(t,x;\cd))$, our conclusion follows for $u(\cd)\in\sD(J(t,x;u(\cd))$. \endpf

\ms

\bc{} \sl Let {\rm(H1), (H3)}, and {\rm(H4)} hold. Then
\bel{D_uJ}D_{u(\cd)}J(t,x;u(\cd))=0,\ee
if and only if
\bel{B=0}\(B(s)+F(s)D(s)\)^\top\wt Y(s)+D(s)^\top\wt Z(s)+S(s)X(s)+R(s)u(s)+r(s)=0.\ee

\ec

\it Proof. \rm It is clear from \rf{D_uJ}. \endpf

\ms

Note that under (H4), the map $u(\cd)\mapsto J(t,x;u(\cd))$ is convex, i.e., the corresponding condition \rf{F_2>0} holds. Therefore, we have the following result.

\bt{} \sl Let {\rm(H1), (H3)} and {\rm(H4)} hold. Then Problem {\rm(LQR)} is open-loop solvable if and only if the following FBSDE is solvable
\bel{FBSDE}\left\{\2n\ba{ll}
\ds dX(s)=[A(s)X(s)+B(s)u(s)+b(s)]ds+[C(s)X(s)+B(s)u(s)+\si(s)]dW(s),\\
\ns\ds d\wt Y(s)=-\[\(A(s)+E(s)I+F(s)C(s)\)^\top\wt Y(s)+\(C(s)+F(s)I\)^\top\wt Z(s)\\
\ns\ds\qq\qq\qq\qq+Q(s)X(s)+S(s)^\top u(s)+q(s)\]ds+\wt Z    (s)dW(s),\qq s\in[t,T],\\
\ns\ds X(t)=x,\q\wt Y(T)=GX(T)+g,\ea\right.\ee
with
\bel{stationary*}\(B(s)+F(s)D(s)\)^\top\wt Y(s)+D(s)^\top\wt Z(s)+S(s)X(s)+R(s)u(s)+r(s)=0,\qq s\in[t,T].\ee

\et

The above result is reduced to the classical one (in \cite{Sun-Yong 2020}) if $E(\cd)=F(\cd)=0$. We see that due to the existence of $E(\cd)$ and $F(\cd)$, our above result is quite different from the classic one. Clearly, FBSDE \rf{FBSDE} is coupled, whose coupling is coming from the {\it stationary condition} \rf{stationary*}. Condition \rf{stationary*} is actually the maximum condition in the maximum principle of optimal control theory.

\ms

\section{The Riccati Differential Equation and the Closed-Loop Solvability}

The goal of this section contains two parts. It will give a characterization for the closed-loop solvability of Problem (LQR).

\subsection{Necessary conditions for Closed-loop solvability}

In this subsection, we are going to establish some necessary conditions for the closed-loop solvability of Problem (LQR). We have the following result.

\bt{necessary} \sl Let {\rm(H1), (H3)--(H4)} hold. Let Problem {\rm(LQR)} be closed-loop solvable a $t\in[0,T]$ with optimal strategy $(\bar\Psi(\cd),\bar v(\cd))\in\BBPsi[t,T]\times\sU^2[t,T]$. Then (suppressing $s$)
\bel{Psi}\left\{\2n\ba{ll}
\ns\ds\bar\Psi=-(R+D^\top PD)^\dag[(B+FD)^\top P+D^\top PC+S],\\
\ns\ds\bar v=-(R+D^\top PD)^\dag[(B+FD)^\top\eta+D^\top P\si+r],\ea\right.\q s\in[t,T],\ee
where $P(\cd)$ is the solution to the solution of the following Riccati differential equation:
\bel{Riccati}\left\{\2n\ba{ll}
\ds\dot P+PA+A^\top P+C^\top PC+Q+EP+F(C^\top P+PC)\\
\ns\ds\qq-[P(B+FD)+C^\top PD+S^\top](R+D^\top PD)^\dag[(B+FD)^\top P+D^\top PC+S]=0,\q s\in[t,T],\\
\ns\ds P(T)=G,\ea\right.\ee
and $\eta(\cd)$ is the solution to the terminal value problem:
\bel{eta}\left\{\2n\ba{ll}
\ns\ds\dot\eta+\{(A+EI+FC)^\top-[P(B+FD)+C^\top PD+S^\top](R+D^\top PD)^\dag(B+FD)^\top\}\eta\\
\ns\ds\qq\qq+\{(C+FI)^\top-[P(B+FD)+C^\top PD+S^\top](R+D^\top PD)^\dag D^\top\} P\si\\
\ns\ds\qq\qq-[P(B+FD)+C^\top PD+S^\top](R+D^\top PD)^\dag r+Pb+q=0,\qq s\in[t,T],\\
\ns\ds\eta(T)=g,\ea\right.\ee
with range conditions
\bel{sR}\left\{\2n\ba{ll}
\ns\ds\sR\((B+FD)^\top P+D^\top PC+S\)\subseteq\sR\(R+D^\top PD\),\\
\ns\ds(B+FD)^\top\eta+D^\top P\si+r\in\sR\(R+D^\top PD\),\ea\right.\ee

\et

\ms

\it Proof. \rm Let Problem (LQR) be closed-loop solvable at $t\in[0,T]$ with $(\bar\Psi(\cd),\bar v(\cd))\in\BBPsi[t,T]\times\sU^2[t,T]$ being an optimal strategy, i.e., \rf{J=inf*} holds. Now, the outcome of this strategy is
\bel{outcome-bar}\bar u(s)=\bar\Psi(s)\bar X(s)+\bar v(s),\qq s\in[t,T].]\ee
It is an open-loop optimal control for the given $x\in\dbR^n$. Thus, the pair $(\bar X(\cd),\bar u(\cd))$, together with some $(\wt Y(\cd),\wt Z(\cd))$ solving the following FBSDE:
\bel{Y_1**}\left\{\2n\ba{ll}
\ds d\bar X(s)=\[A(s)\bar X(s)+B(s)\(\bar\Psi(s)\bar X(s)+\bar v(s)\)+b(s)\]ds\\
\ns\ds\qq\qq\qq+\[C(s)\bar X(s)+D(s)\(\Psi(s)\bar X(s)+\bar v(s)\)+\si(s)\]dW(s),\q s\in[t,T],\\
\ns\ds d\wt Y(s)=-\[\(A(s)+E(s)I+F(s)C(s)\)^\top\wt Y(s)+\(C(s)+F(s)I\)^\top\wt Z(s)\\
\ns\ds\qq\qq\qq+Q(s)\bar X(s)+S(s)^\top\(\Psi(s)\bar X(s)+\bar v(s)\)+q(s)\]ds+\wt Z(s)dW(s),\q s\in[t,T],\\
\ns\ds\bar X(t)=x,\q\wt Y(T)=G\bar X(T)+g,\ea\right.\ee
with the stationary condition:
\bel{stationary}\(B(s)+F(s)D(s)\)^\top\wt Y(s)+D(s)^\top\wt Z(s)+S(s)\bar X(s)+R(s)\(\Psi(s)\bar X(s)+\bar v(s)\)+r(s)=0,\q s\in[t,T].\ee
We can write \rf{Y_1**} as
\bel{wt Y}\left\{\2n\ba{ll}
\ds d\bar X(s)=\[\(A(s)+B(s)\bar\Psi(s)\)\bar X(s)+B(s)\bar v(s)+b(s)\]ds\\
\ns\ns\qq\qq\qq+\[\(C(s)+D(s)\bar\Psi(s)\)\bar X(s)+D(s)\bar v(s)+\si(s)\]dW(s),\\
\ds d\wt Y(s)=-\[\(A(s)+E(s)I+F(s)C(s)\)^\top\wt Y(s)+\(C(s)+F(s)I\)^\top\wt Z(s)\\
\ns\ds\qq\qq\qq+\(Q(s)+S(s)^\top\bar\Psi(s)\)\bar X(s)+S(s)^\top\bar v(s)+q(s)\]ds+\wt Z(s)dW(s),\qq s\in[t,T],\\
\ns\ds\bar X(t)=x,\q\wt Y(T)=G\bar X(T)+g,\ea\right.\ee
with the stationary condition
\bel{stationary-wt}\(B(s)+F(s)D(s)\)^\top\wt Y(s)+D(s)^\top\wt Z(s)+\(S(s)+R(s)\bar\Psi(s)\)\bar X(s)+R(s)\bar v(s)+r(s)=0.\ee
Similarly, under the same affine strategy, the outcome, denoted by $\bar u_0(\cd)$, is an open-loop optimal control for the initial state $0$. Thus, the following holds: (with the optimal pair denoted by $(\bar X_0(\cd),\bar\Psi(\cd)\bar X_0(\cd)+\bar v(\cd))$, and the adapted solution $(\wt Y_0(\cd),\wt Z_0(\cd))$ of the BSDE)
\bel{Y*}\left\{\2n\ba{ll}
\ds d\bar X_0(s)=\[\(A(s)+B(s)\bar\Psi(s)\)\bar X_0(s)+B(s)\bar v(s)+b(s)\]ds\\
\ns\ns\qq\qq\qq+\[\(C(s)+D(s)\bar\Psi(s)\)\bar X_0(s)+D(s)\bar v(s)+\si(s)\]dW(s),\\
\ds d\wt Y_0(s)=-\[\(A(s)+E(s)I+F(s)C(s)\)^\top\wt Y_0(s)+\(C(s)+F(s)I\)^\top\wt Z_0(s)\\
\ns\ds\qq\qq\qq+\(Q(s)+S(s)^\top\bar\Psi(s)\)\bar X_0(s)+S(s)^\top\bar v(s)+q(s)\]ds+\wt Z_0(s)dW(s),\q s\in[t,T],\\
\ns\ds\bar X_0(t)=0,\q\wt Y_0(T)=G\bar X_0(T)+g,\ea\right.\ee
with the corresponding stationary condition:
\bel{stationary***}\(B(s)+F(s)D(s)\)^\top\wt Y_0(s)+D(s)^\top\wt Z_0(s)+\(S(s)+R(s)\bar\Psi(s)\)\bar X_0(s)+R(s)\bar v(s)+r(s)=0.\ee
Let
\bel{h X}\h X(s)=\bar X(s)-\bar X_0(s).\ee
Then we obtain the following homogeneous system
\bel{Y*}\left\{\2n\ba{ll}
\ds d\h X(s)=\[\(A(s)+B(s)\bar\Psi(s)\)\h X(s)\]ds+\[\(C(s)+D(s)\bar\Psi(s)\)\h X(s))\]dW(s),\\
\ds d\h Y(s)=-\[\(A(s)+E(s)I+F(s)C(s)\)^\top\h Y(s)+\(C(s)+F(s)I\)^\top\h Z(s)\\
\ns\ds\qq\qq\qq+\(Q(s)+S(s)^\top\bar\Psi(s)\)\h X(s)\]ds+\h Z(s)dW(s),\q s\in[t,T],\\
\ns\ds\h X(t)=x,\q\h Y(T)=G\h X(T),\ea\right.\ee
with the difference of the stationary conditions:
\bel{diff}\(F(s)D(s)^\top+B(s)^\top\)\h Y(s)+D(s)^\top\h Z(s)+\(S(s)+R(s)\bar\Psi(s)\)\h X(s)=0.\ee
Clearly, \rf{Y*}--\rf{diff} is a homogeneous FBSDE. It has a unique adapted solution $(\h X(\cd),\h Y(\cd),\h Z(\cd))$ (for any $x\in\dbR^n$). We may let
matrix-valued process FBSDE as follows, which is solvable:
\bel{XYZ}\left\{\2n\ba{ll}
\ds d\dbX(s)=\[\(A(s)+B(s)\bar\Psi(s)\)\dbX(s)\]ds+\[\(C(s)+D(s)\bar\Psi(s)\)\dbX(s))\]
dW(s),\\ [2mm]
\ds d\dbY(s)=-\[\(A(s)+E(s)I+F(s)C(s)\)^\top\dbY(s)+\(C(s)+F(s)I\)^\top\dbZ(s)\\
\ns\ds\qq\qq+\(Q(s)+S(s)^\top\bar\Psi(s)\)\dbX(s)\]ds+\dbZ(s)dW(s),\q s\in[t,T],\\
\ns\ds\dbX(t)=I,\q\dbY(T)=G\dbX(T),\ea\right.\ee
with the stationary condition:
\bel{stationary}\(F(s)D(s)^\top+B(s)^\top\)\dbY(s)+D(s)^\top\dbZ(s)
+\(S(s)+R(s)\bar\Psi(s)\)\dbX(s)=0.\ee
Note
$$\ba{ll}
\ns\ds0=d(\dbX\dbX^{-1})=[(A+B\bar\Psi)\dbX ds+(C+D\bar\Psi)\dbX dW]\dbX^{-1}+\dbX(\a ds+\b dW)+(C+D\bar\Psi)\dbX\b ds\\
\ns\ds\q=[(A+B\bar\Psi)+\dbX\a+(C+D\bar\Psi)\dbX\b]ds+[(C+D\bar\Psi)+\dbX\b] dW.\ea$$
Then
$$\b=-\dbX(C+D\bar\Psi),\qq\a=\dbX^{-1}[(C+D\bar\Psi)^2-(A+B\bar\Psi)].$$
Hence, $\dbX^{-1}$ exists and it satisfies
\bel{X^{-1}}\left\{\2n\ba{ll}
\ns\ds d[\dbX^{-1}]=\dbX^{-1}[(C+D\bar\Psi)^2-(A+B\bar\Psi)]ds-\dbX^{-1}(C+D\bar\Psi)dW,\\
\ns\ds\dbX^{-1}(t)=I.\ea\right.\ee
Define
\bel{P}P=\dbY\dbX^{-1},\qq\L=\dbZ\dbX^{-1}.\ee
By \rf{stationary}, we have
\bel{sttionary2}\(F(s)D(s)^\top+B(s)^\top\)P(s)+D(s)^\top\L(s)
+\(S(s)+R(s)\bar\Psi(s)\)=0.\ee
Then, using the It\^o's formula,
\bel{dP}\ba{ll}
\ns\ds dP=d[\dbY\dbX^{-1}]=-[(A+EI+FC)^\top P+(C+FI)^\top\L+(Q+S^\top\bar\Psi)]ds+\L dW\\
\ns\ds\qq\qq+P[(C+D\bar\Psi)^2-(A+B\bar\Psi)]ds-P(C+D\bar\Psi)dW
-\L(C+D\bar\Psi)ds\\
\ns\ds\q~=-[(A+EI+FC)^\top P+(C+FI)^\top\L+(Q+S^\top\bar\Psi)
+P(A+B\bar\Psi)\\
\ns\ds\qq\qq-P(C+D\bar\Psi)^2+\L(C+D\bar\Psi)]ds+[\L-P(C+D\bar\Psi)]dW\ea\ee
Since $P(\cd)$ is deterministic, we must have (the diffusion vanishes)
\bel{L}\L=P(C+D\bar\Psi).\ee
Then \rf{dP} can be written as
$$\ba{ll}
\ns\ds dP=d[\dbY\dbX^{-1}]=-[(A+EI+FC)^\top P+(C+FI)^\top P(C+D\bar\Psi)+(Q+S^\top\bar\Psi)
+P(A+B\bar\Psi)]ds\\
\ns\ds\q~=-[PA+A^\top P+C^\top PC+Q+(EI+FC)^\top P+FPC+C^\top PD\bar\Psi+FPD\bar\Psi+S^\top\bar\Psi
+PB\bar\Psi]ds\\
\ns\ds\q~=-[\BQ(P)+(EI+FC)^\top P+FPC+\BS(P)^\top\bar\Psi+FPD\bar\Psi]ds\\
\ns\ds\q~=-\{\BQ(P)+EP+F(C^\top P+PC)+[\BS(P)^\top+FPD]\bar\Psi\}ds,\ea$$
where
\bel{BQ}\left\{\2n\ba{ll}
\ns\ds\BQ(s,P)=PA(s)+A(s)^\top P+C(s)^\top PC(s)+Q(s),\\
\ns\ds\BS(s,P)=B(s)^\top P+D(s)^\top PC(s)+S(s),\ea\right.\qq (s,P)\in[t,T]\times\dbS^n.\ee
In what follows, we also let
\bel{BR}\BR(s,P)=R(s)+D(s)^\top PD(s),\qq(s,P)\in[t,T]\times\dbS^n.\ee
Since (from the homogeneous stationary condition \rf{stationary})
\bel{sttionary3}0=(FD^\top P+B^\top P+D^\top P(C+D\bar\Psi)+
(S+R\bar\Psi)=FD^\top P+\BS(P)+\BR(P)\bar\Psi,\ee
it implies
\bel{range*}\sR\big(FD^\top P+\BS(P)\big)\subseteq\sR\big(\BR(P)\big),\ee
Thus
\bel{bar Psi}\bar\Psi\equiv\bar\Psi(P)=-\BR(P)^\dag[\BS(P)+FD^\top P].\ee
Consequently, the following Riccati differential equation is solvable:
\bel{Riccati*}\left\{\2n\ba{ll}
\ds\dot P+\BQ(P)+EP+F(C^\top P+PC)-[\BS(P)^\top+FPD]\BR(P)^\dag[\BS(P)+FD^\top P]=0,\\
\ns\ds P(T)=G,\ea\right.\ee
which is the same as \rf{Riccati}. Now, we set
\bel{wteta}\left\{\2n\ba{ll}
\ns\ds\wt\eta(s)=\wt Y(s)-P(s)\bar X(s),\\
\ns\ds\wt\z(s)=\wt Z(s)-P(s)[(C(s)+D(s)\bar\Psi(s))\bar X(s)+D(s)\bar v(s)+\si(s)].\ea\right.\ee
Then (note \rf{wt Y})
$$\ba{ll}
\ns\ds d\wt\eta=d\wt Y-[\dot PX]ds-PdX\\
\ns\ds\q=-[(A+EI+FC)^\top\wt Y+(C+FI)^\top\wt Z+(Q+S^\top\bar\Psi)\bar X+S^\top\bar v+q]ds+\wt ZdW(s)\\
\ns\ds\qq\qq-\dot P\bar Xds-P\{[(A+B\bar\Psi)\bar X+B\bar v+b]ds+[(C+D\bar\Psi)\bar X+D\bar v+\si]dW(s)\}\\
\ns\ds\q=-\{\dot P\bar X+(A+EI+FC)^\top\wt Y+(C+FI)^\top\wt Z+(Q+S^\top\bar\Psi)\bar X+S^\top\bar v+q\\
\ns\ds\qq\qq+P[(A+B\bar\Psi)\bar X+B\bar v+b]\}ds+\{\wt Z-P[(C+D\bar\Psi)\bar X+D\bar v+\si]\}dW(s)\\
\ns\ds\q=-\{\dot P\bar X+(A+EI+FC)^\top(P\bar X+\eta)+(C+FI)^\top\{\wt\z+P[(C+D\bar\Psi)\bar X+D\bar v+\si]\}\\
\ns\ds\qq\qq+(Q+S^\top\bar\Psi)\bar X+S^\top\bar v+q+P[(A+B\bar\Psi)\bar X+B\bar v+b]\}ds+\z dW(s)\\
\ns\ds\q=-\{[\dot P+(A+EI+FC)^\top P+(C+FI)^\top PC+PA+Q]\bar X+(A+EI+FC)^\top\wt\eta\\
\ns\ds\qq\qq+[PB+(C+FI)^\top PD+S^\top]\bar\Psi\bar X+(C+FI)^\top\wt\z+[PB+(C+FI)^\top PD+S^\top)\bar v+b]\}\\
\ns\ds\qq\qq+(C+FI)^\top P\si+q+Pb\}ds+\wt\z dW(s)\\
\ns\ds\q=-\{[\dot P+EP+F(C^\top P+PC)+\BQ(P)+(\BS(P)^\top+FPD)\bar\Psi]\bar X+(A+EI+FC)^\top\wt\eta\\
\ns\ds\qq\qq+(C+FI)^\top\wt\z+(FPD+\BS(P)^\top)\bar v+b]+(C+FI)^\top P\si+q+Pb\}ds+\wt\z dW(s)\\
\ns\ds\q=-\{(A+EI+FC)^\top\wt\eta+(C+FI)^\top\wt\z+(FPD+\BS(P)^\top)\bar v+b]+(C+FI)^\top P\si+q+Pb\}ds\\
\ns\ds\qq\qq\qq\qq\qq\qq\qq\qq\qq\qq\qq\qq\qq\qq+\wt\z dW(s)\ea$$
From the stationary condition \rf{stationary}, we have (note \rf{wteta})
$$\ba{ll}
\ns\ds0=(B+FD)^\top\wt Y+D^\top\wt Z+(S+R\bar\Psi)\bar X+R\bar v+r\\
\ns\ds\q=(B+FD)^\top(P\bar X+\wt\eta)+D^\top[PC\bar X+PD(\bar\Psi\bar X+\bar v)+P\si+\wt\z]+(S+R\bar\Psi)\bar X+R\bar v+r\\
\ns\ds\q=(B+FD)^\top P\bar X+(B+FD)^\top\wt\eta+D^\top PC\bar X+D^\top PD\bar\Psi\bar X+D^\top PD\bar v+D^\top P\si+D^\top\wt\z\\
\ns\ds\qq\qq+S\bar X+R\bar\Psi\bar X+R\bar v+r\\
\ns\ds\q=[\BS(P)+FD^\top P+\BR(P)\bar\Psi]\bar X+(B+FD)^\top\wt\eta+\BR(P)\bar v+D^\top P\si+D^\top\wt\z+r\\
\ns\ds\q=\BR(P)\bar v+(B+FD)^\top\wt\eta+D^\top P\si+D^\top\wt\z+r.\ea$$
This implies
\bel{r in R}(B+FD)^\top\wt\eta+D^\top P\si+D^\top\wt\z+r\in\sR\big(\BR(P)\big),\ee
and
\bel{v}\bar v=-\BR(P)^\dag\[(B+FD)^\top\wt\eta+D^\top P\si+D^\top\wt\z+r\].\ee
Hence, $(\wt\eta(\cd),\wt\z(\cd))$ is the adapted solution of BSDE:
\bel{eta*}\left\{\2n\ba{ll}
\ds d\wt\eta=-\g ds+\wt\z dW(s),\q s\in[t,T],\\
\ns\ds\wt\eta(T)=g,\ea\right.\ee
with
\bel{g}\ba{ll}
\ns\ds\g=[(A+EI+FC)^\top-(\BS(P)^\top+FPD)\BR(P)^\dag[(B+FD)^\top]\wt\eta\\
\ns\ds\qq\qq+\{(C+FI)^\top-[\BS(P)^\top+FPD]\BR(P)^\dag D^\top\}\wt\z+Pb+q\\
\ns\ds\qq\qq+\{(C+FI)^\top-[\BS(P)^\top+FPD]\BR(P)^\dag D^\top\}P\si-[\BS(P)^\top+FPD]\BR(P)^\dag r.\ea\ee
Now, since \rf{eta} is a BSDE with deterministic coefficients and deterministic  terminal state, we must have $\wt\z=0$. Thus,
\rf{g} reads
\bel{g*}\ba{ll}
\ns\ds\g=[(A+EI+FC)^\top-(\BS(P)^\top+FPD)\BR(P)^\dag[(B+FD)^\top]\wt\eta\\
\ns\ds\qq+\{(C+FI)^\top-[\BS(P)^\top+FPD]\BR(P)^\dag D^\top\}P\si-[\BS(P)^\top+FPD]\BR(P)^\dag r+Pb+q.\ea\ee
Then, combining the above with the uniqueness, we obtain ODE \rf{eta}. \endpf

\subsection{Sufficient conditions for closed-loop solvability}

In this subsection, we are going to give a sufficient condition for the closed-loop solvability of Problem (LQR). The ideas come from the so-called the Four Step-Scheme found in \cite{Ma-Protter-Yong 1994, Ma-Yong 1999} for FBSDEs. Our main result is the following.

\bt{} \sl Let {\rm(H1), (H3)--(H4)} hold. Let the Riccati differential equation \rf{Riccati} (or \rf{Riccati*}) admit a unique solution $P(\cd)$, and let $\eta(\cd)$ be the solution of the terminal value problem \rf{eta} (or
\rf{eta*}--\rf{g}), with range conditions \rf{sR}. Then Problem {\rm(LQR)} is closed-loop solvable at any $t\in[0,T]$ with the optimal closed-loop strategy given by \rf{Psi}.
\et

\it Proof. \rm Let $(\bar\Psi(\cd),\bar v(\cd))\in\BBPsi[t,T]\times\sU^2[t,T]$ be any given strategy. For any initial state $x\in\dbR^n$, let
$$\bar u(\cd)=\bar\Psi(\cd)\bar X(\cd)+\bar v(\cd)$$
be a corresponding outcome of this strategy (with $\bar X(\cd)$ being the corresponding state process). Also, let
\bel{wt Y**}\left\{\2n\ba{ll}
\ns\ds\bar Y(s)=P(s)\bar X(s)+\eta(s),\\
\ns\ds\bar Z(s)=P(s)\[C(s)\bar X(s)+D(s)\bar u(s)+\si(s)\],\ea
\right.\qq s\in[t,T],\ee
with $P(\cd)$ being the solution to the Riccati dyifferential equation \rf{Riccati} and $\eta(\cd)$ being the solution of terminal value problem \rf{eta}.
Then, we claim that $(\bar X(\cd),\bar u(\cd),\bar Y(\cd),\bar Z(\cd))$ is the adapted solution of the following FBSDE:
\bel{FBSDE}\left\{\2n\ba{ll}
\ds d\bar X(s)=\[A(s)\bar X(s)+B(s)\bar u((s)+b(s)\]ds+\[C(s)\bar X(s)+D(s)\bar u(s)+\si(s)\]dW(s),\\
\ns\ds d\bar Y(s)=-\[\(A(s)+E(s)I+F(s)C(s)\)^\top\bar Y(s)+\(C(s)+F(s)I\)^\top\bar Z(s)\\
\ns\ds\qq\qq\qq\qq+Q(s)\bar X(s)+S(s)^\top\bar u(s)+q(s)\]ds+\bar Z    (s)dW(s),\qq s\in[t,T],\\
\ns\ds\bar X(t)=x,\q\bar Y(T)=G\bar X(T)+g,\ea\right.\ee
with the stationary condition
\bel{stationary*}[B(s)+F(s)D(s)]^\top\bar Y(s)+D(s)^\top\bar Z(s)+S(s)\bar X(s)+R(s)\bar u(s)+r(s)=0.\ee
To solve \rf{FBSDE}--\rf{stationary*} (by means of Four Step Scheme \cite{Ma-Protter-Yong 1994,Ma-Yong 1999}), we set
$$\bar Y(s)=\bar P(s)\bar X(s)+\bar\eta(s),\qq s\in[t,T].$$
Then
$$\bar Y(T)=G\bar X(T)+g.$$
Next, we need
$$\ba{ll}
\ns\ds-[(A+EI+FC)^\top\bar Y+(C+FI)^\top\bar Z+Q\bar X+S^\top\bar u+q]ds+\bar ZdW\\
\ns\ds=d\bar Y=\dot P\bar Xds+P\[\(A\bar X+B\bar u+b\)ds+\(C\bar X+D\bar u+\si\)dW\]+\dot\eta ds.\ea$$
Thus, the requirement for $\bar Z(s)$ in \rf{wt Y**} holds. Next, we require the drifts on the both sizes of \rf{wt Y**} are equal:
$$\ba{ll}
\ns\ds0=\dot PX+P(A\bar X+B\bar u+b)+\dot\eta+(A+EI+FC)^\top\bar Y+(C+FI)^\top\bar Z+Q\bar X+S^\top\bar u+q\\
\ns\ds\q=\dot PX+P(A\bar X+B\bar u+b)+\dot\eta+(A+EI+FC)^\top(P\bar X+\eta)\\
\ns\ds\qq+(C+FI)^\top P(C\bar X+D\bar u+\si)+Q\bar X+S^\top\bar u+q\\
\ns\ds\q=[\dot P+PA+A^\top P+C^\top PC+Q+EP+F(C^\top P+PC)]\bar X\\
\ns\ds\qq+(PB+C^\top PD+S^\top+FPD)\bar u+Pb+\dot\eta+(A+EI+FC)^\top\eta+(C+FI)^\top P\si+q\\
\ns\ds\q=[\dot P+\BQ(P)+EP+F(C^\top P+PC)]\bar X+[\BS(P)^\top+FPD]\bar u\\
\ns\ds\qq+(A+EI+FC)^\top\eta+(C+FI)^\top P\si+Pb+q+\dot\eta\ea$$
Note we also need:
$$\ba{ll}
\ns\ds0=(B+FD)^\top\bar Y+D^\top\bar Z+S\bar X+R\bar u+r\\
\ns\ds\q=(B+DF)^\top(P\bar X+\eta)+D^\top P(C\bar X+D\bar u+\si)+S\bar X+R\bar u+r\\
\ns\ds\q=(B^\top P+D^\top PC+S+FD^\top P)\bar X+(R+D^\top PD)\bar u+(B+FD)^\top\eta+D^\top P\si+r\\
\ns\ds\q=[\BS(P)+FD^\top P]\bar X+\BR(P)\bar u+(B+FD)^\top\eta+D^\top P\si+r.\ea$$
Under the range conditions \rf{sR}, we have
$$\bar u=-\BR(P)^\dag[\BS(P)+FD^\top P]\bar X-\BR(P)^\dag[(B+FD)^\top\eta+D^\top P\si+r].$$
With such a control, we have
$$\ba{ll}
\ns\ds0=[\dot P+\BQ+EP+F(CP^\top+PC)]\bar X-[\BS(P)^\top+FPD]\BR(P)^\dag[\BS(P)+FD^\top P]\bar X\\
\ns\ds\qq-\BR(P)^\dag[\BS(P)^\top+FPD][(B+FD)^\top\eta+D^\top P\si+r]\\
\ns\ds\qq+(A+EI+FC)^\top\eta+(C+FI)^\top P\si+Pb+q+\dot\eta\\
\ns\ds\q=\{(A+EI+FC)^\top-\BR(P)^\dag[\BS(P)^\top+FPD][(B+FD)^\top\}\eta\\
\ns\ds\qq+\{(C+FI)^\top-\BR(P)^\dag[\BS(P)^\top+FPD]D^\top\}P\si+D^\top P\si+r]\\
\ns\ds\qq-\BR(P)^\dag r+Pb+q+\dot\eta.\ea$$
Thus, \rf{FBSDE}--\rf{stationary*} holds. This shows that the FBSDE
\rf{wt Y} is decoupled, for any $x\in\dbR^n$. Then, we obtain the open-loop solvability of Problem (LQR) for any $x\in\dbR^n$, under the same strategy $(\Psi(\cd),v(\cd))\in\BBPsi[t,T]\times\sU^2[t,T]$ given by \rf{Psi}.
Hence, Problem (LQR) is closed-loop solvable at $t\in[0,T]$. \endpf

\section{Conclusions}

We formulated a stochastic LQ problem with recurve cost functional. It is naturally involved BSDE in $L^1$, which leads to some technicality in the correct formulation of Problem (LQ). We defined the open-loop and closed-loop solvability and their characteristics. It is proved that for Problem (LQR), its open-loop solvable at $(t,x)$ if and only if the corresponding FBSDE (together with the stationary condition) is solvable on $[t,T]$; and its closed-loop solvable at $t$ if and only if the Riccati differential equation admits a solution, and a terminal value problem is solvable.

\ms

We expect that by assuming uniform convexity conditions, Problem (LQR) is both open-loop and closed-loop solvable. We are not able to prove that for the time being. Thus, it is left open.

\end{document}